\documentclass[graybox,12pt]{svmult}

\usepackage{mathptmx}       
\usepackage{helvet}         
\usepackage{courier}        

\usepackage{makeidx}         
\usepackage{graphicx}        
\usepackage{multicol}        
\usepackage{dsfont}
\usepackage{amsmath}
\usepackage{amsfonts}
\usepackage{color}     
\usepackage{url}
\usepackage{verbatim}
\usepackage{enumerate}

\makeindex             

\usepackage{amssymb}
\usepackage{amsmath}
\usepackage{epsfig}
\usepackage{mathptmx}       
\usepackage{helvet}         
\usepackage{courier}        
\usepackage{type1cm}        

\usepackage{makeidx}         
\usepackage{graphicx}        
\usepackage{multicol}        
\usepackage[bottom]{footmisc}

\usepackage{amsbsy}
\usepackage{MnSymbol}
\usepackage{amsfonts}
\usepackage{amsopn}
\usepackage{dsfont}
\usepackage[english]{babel}
\usepackage[applemac]{inputenc}
\usepackage[colorlinks=true,breaklinks=true]{hyperref}
\usepackage[T1]{fontenc}
\usepackage{xcolor}
\usepackage{enumerate}
\usepackage{latexsym,amsmath,eufrak}
\usepackage{xcolor}
\usepackage{stackrel}
\newcommand{\Rset}{\mathds{R}}
\newcommand{\Cset}{\mathds{C}}
\newcommand{\Hset}{\mathds{H}}

\newcommand{\oi}{(0,\infty)}
\newcommand{\CF}{\mathcal{CF}}
\newcommand{\CBF}{\mathcal{CBF}}
\newcommand{\CBe}{\CBF_{\!\!e}}
\newcommand{\TBF}{\mathcal{TBF}}
\newcommand{\ST}{\mathcal{S}}
\newcommand{\BF}{\mathcal{BF}}
\newcommand{\Nev}{\mathcal{N}}
\newcommand{\Four}{\mathfrak{F}}

\newcommand{\Lle}{\mathcal{LE}}

\newcommand{\simdis}{\stackrel{d}{=}}
\newcommand{\Ebf}{\mbox{\bf e}}
\newcommand{\Ubf}{\mbox{\bf u}}
\newcommand{\er}{\mathbb{E}}
\newcommand{\pr}{\mathbb{P}}

\newcommand{\II}{{\rm  1~\hspace{-1.1ex}l}}
\newcommand{\qqr}{\mathrm{q}}
\newcommand{\aar}{\mathrm{a}}
\newcommand{\bbr}{\mathrm{b}}
\newcommand{\ccr}{\mathrm{c}}
\newcommand{\ddr}{\mathrm{d}}

\newcommand{\Xs}{{\scriptscriptstyle X}}
\newcommand{\Xts}{{\scriptscriptstyle \tilde{X}}}
\newcommand{\Cs}{{\scriptscriptstyle C}}
\newcommand{\Ss}{{\scriptscriptstyle S}}
\newcommand{\Ts}{{\scriptscriptstyle T}}
\newcommand{\Cts}{{\scriptscriptstyle \tilde{C}}}
\newcommand{\Sts}{{\scriptscriptstyle \tilde{S}}}
\newcommand{\Tts}{{\scriptscriptstyle \tilde{T}}}
\newcommand{\as}{{\scriptscriptstyle \alpha}}
\newcommand{\Ls}{\scalebox{0.5}{$\mbox{\bf  L}$}}

\newcommand{\Lap}{\mathfrak{L}}
\newcommand{\Fou}{\mathfrak{F}}

\newcommand{\Ei}{{\rm Ei}}
\newcommand{\Li}{{\rm Li}}
\newcommand{\li}{{\mathfrak{Li}}}
\newcommand{\ei}{{\mathfrak{Ei}}}

\newcommand\rsmraise[1]{%
  \ifx#1\displaystyle .8\else
    \ifx#1\textstyle .8\else
      \ifx#1\scriptstyle .6\elseis
        .45%
      \fi
    \fi
  \fi}

\newenvironment{proofs}{\noindent {\bf Proof }}{\ \ \ $\square$\\}

\begin{document}
\title*{Integral transforms related to Nevanlinna-Pick functions  from an analytic, probabilistic and free-probability point of view}
\author{Wissem Jedidi  \thanks{Department of Statistics \& OR, King Saud University, P.O. Box 2455, Riyadh 11451, Saudi Arabia, \email{wjedidi@ksu.edu.sa}, \email{nuha.altaymani@gmail.com}} \thanks{Universit\'e de Tunis El Manar, Facult\'e des Sciences de Tunis, D\'epartement de Math\'ematiques, Laboratoire d'Analyse Math\'ematiques et Applications LR11ES11.  2092 - El Manar I, Tunis, Tunisia, \email{wissem.jedidi@fst.utm.tn}}, Zbigniew J.  Jurek \thanks{University of Wroc\l aw,  Poland, \email{zjjurek@math.uni.wroc.pl}} and Nuha Taymani $^{\ast}$}
\maketitle
\date{\today}
\abstract{We establish a new connection between the class of Nevanlinna-Pick functions and the one of the exponents associated to spectrally negative L\'evy processes.
As a consequence, we compute the characteristics related to some hyperbolic functions and we show a property of temporal  complete monotonicity, similar to the one obtained  via the Lamperti transformation by Bertoin \& Yor  ({\it On subordinators, self-similar Markov processes and some factorizations of the exponential variable}, Elect. Comm. in Probab., vol. 6, pp. 95--106, 2001) for self-similar Markov processes. More precisely, we show the remarkable fact that for a subordinator  $\xi$, the function $t \mapsto t^n \, \er[\xi_t^{-p}]$ is , depending on the values of the exponents $n=0,1,2,\; p>-1$, or a Bernstein function or a  completely monotone function. In particular, $\xi$ is the inverse  time  subordinator  of a spectrally negative L\'evy process, if, and only if, for some $\,p\geq 1$, the function $t \mapsto t \, \er[\xi_t^{-p}]$  is a  Stieltjes transform. Finally, we clarify to which extent Nevanlinna-Pick functions are related to free-probability and to Voiculescu transforms, and we provide an inversion procedure.
}
\section{Introduction}
In complex analysis, analytic functions that preserve the upper complex half-plane play an important role. According to Nevanlinna-Pick  theorem , these are precisely functions $F$ that admit the (unique) representation:
\begin{equation}\label{ffz}
F(z)= \aar+\bbr\; z+\int_{\Rset}\frac{zx-1}{z+x}\rho(dx), \quad \mbox{for} \; z \in \Cset \setminus \Rset, \; \mbox{where} \;   \aar\in \Rset,\, \bbr\ge 0, \; \mbox{and $\rho$ is a finite Borel measure.}
\end{equation}
We will denote by $\Nev$  the set of all {\it Nevanlinna-Pick functions} while $\Nev(A)$, for a Borel set $A$, stands for
$$\Nev(A):=\{\,F \in \Nev \;\mbox{ such that in  \eqref{n1}}, \;support(\rho)\subset A \}.$$

On the other hand, there are {\it complete Bernstein functions}, cf. Appendix, i.e. those functions  $f:(0, \infty)\to \Rset$ that analytically extend to $z \in \Cset \setminus (-\infty,0]$ and have a form as Nevanlinna-Pick functions as in \eqref{ffz}. Furthermore, in the free-probability theory (cf. Appendix) an infinitely divisible
probability measure $\nu$ has its {\it Voiculescu transform} $V_{\nu}(z)$ as follows
\begin{equation}\label{vovo}
V_\nu(z)= \aar + \int_{\Rset}\frac{zx-1}{z+x}m(dx), \quad \mbox{for} \; z \in \Cset \setminus \Rset,
\end{equation}
for a uniquely determined finite measure $m$. In this paper, we will not use the analyticity property of Nevanlinna-Pick functions or Voiculesu transforms. In fact, from the values of $F(it)$ or $V_{\nu}(it), \;t \in \Rset$, we can retrieve the representing parameters, in particular, the $\rho$-measures, cf. \eqref{recup} below.\\

Recall that in the (classical) probability we have that a measure $\mu$ is infinitely divisible if and only if  its Fourier transform (or characteristic function) $\mathfrak{F}_{\mu}(t)$ is of a form
\begin{equation}\label{pvovo}
\mathfrak{F}_{\mu}(s): =\int_{\Rset}e^{i\,s\,x}\mu(dx)=\exp\Big[i\,u\,\aar - \bbr^2 u^2 +\int_{\Rset}\left(e^{i\,s\,x}-1-\frac{itx}{1+x^2}\right)\frac{1+x^2}{x^2}m(dx)\Big],
\end{equation}
where $\aar \in \Rset,$ $\bbr\geq 0,$ and $m$ is a finite Borel measure.  The above is the famous {\it L\'evy-Khintchine formula} for infinitely divisible characteristic functions and those are closely related to L\'evy stochastic process. Note that \eqref{vovo}  and \eqref{pvovo} give a bijection between the classical infinite divisibility and the free-probability infinite divisibility. Namely, $\nu$ in \eqref{vovo} is a free-probability analogue of $\mu$ in the classical probability (when $\bbr=0$).\\

In this paper, Theorem \ref{vraiz},  gives some equivalent characterizations of $F\in \Nev(A)$, among them, there is a way to associate $F$ with a Cauchy transform. In Proposition \label{prop-moment} and Theorem \ref{temporalcm},  we establish  a temporal monotonicity property for subordinators. Finally, in Propositions \ref{pro1} and \ref{pro4},  are found  Voiculescu transforms of free-probability analogues of the classical infinitely divisible hyperbolic functions. In particular, the $\rho$-measures in \eqref{vovo}   are given as a composition of  Laplace and Fourier transforms.\\

In order to precise our main results, we need to formalize our objects.
The  right, left,  upper and lower complex half-plane are denoted by
$$\Cset^+:=\{z\in\Cset: \Re(z)>0\},\quad \Cset^-:=-\Cset^+, \quad \Hset^+:=\{z\in\Cset:\Im(z)>0\}, \quad \Hset^-:=-\Hset^+.$$
The Alexandrov compactification of $\Cset\cup \{\infty\}$ is denoted by $\overline{\Cset}$. For given $\alpha,N > 0$ and $M\in \Rset$, we introduce the angular sectors and cones in $\Hset^+$:
\begin{equation}
\Gamma_\alpha(M)  :=   \{z=a+ib \in\Hset^+\;:\;    \alpha b >  |a -M|\},\qquad \Gamma_{\alpha} := \Gamma_\alpha(0)=\Gamma_{\alpha}(M) -M, \qquad
\Gamma_{\alpha,N} := \{z=a+ib\in \Gamma_{\alpha}\;:\; b >N\}.\label{reg}
\end{equation}

In classical  complex analysis, one of  the  fundamental results is the  integral representation of \emph{analytic functions} defined on  $\Hset^+$, or $\Cset {\footnotesize \setminus} \Rset=\Hset^+\cup\Hset^-$,
that preserve $\Hset^+$. If such an analytic function, say $F$, is only defined on $\Hset^+$, then we can extend it onto $\Cset {\footnotesize \setminus} \Rset$ by setting
$F(z)= \overline{F(\overline{z})}$, if $z \in \Hset^-$. Then, $F$ admits the (unique) canonical form: for $z\in  \Cset {\footnotesize \setminus} \Rset$,
\begin{eqnarray}
F(z)&=&\aar + \bbr z+ \int_{\Rset} \frac{zx -1}{z+x}\rho(dx)\label{n1}\\
&=& \aar +\bbr z+ \int_{\Rset}\left(\frac{x}{1+x^2} -\frac{1}{x+z} \right) (1+x^2)\rho(dx).\label{n2}
\label{star}\end{eqnarray}
where $\aar\in \Rset$ is a real number, $\bbr\geq 0$ and $\rho$ is a finite (Borel) measure on $\Rset$.
Note that
\begin{equation}\label{recup}
\aar=\Re\left(F(i)\right), \quad \bbr=\lim_{u\to \infty} F(iu)/(iu)\quad  \mbox{and}\quad \rho(\Rset)= \Im \left(F(i)\right)-\bbr,
\end{equation}
that for every measure $\rho$ we have the inversion formula, cf. \cite[p.126]{Akhi} or Lang \cite[p.380]{Lang}:
\begin{equation}\int_{(u,v]} (1+x^2) \rho(dx)= \lim_{\epsilon \to 0^+}\frac{1}{\pi}\int_{(u,v]} \Im\left(F(-x+ i \epsilon)\right)\, dx,\quad \mbox{whenever $\rho(\{u,v\})=0$},
\label{stars}\end{equation}
and that the half-plane preservation property could be seen from property {\bf(P1)} in Section \ref{free} , of the function $f$ in \eqref{mob} below. Those function $F$ are coined in the literature as  {\it Pick or Nevanlinna-Pick functions}. cf. Akhiezer \cite{Akhi}, Bondesson \cite[p.21]{Bon},  or Schilling, Song and Vondra\v{c}ek \cite[Chapter 6]{SSV}. The (unique) representing measure $\rho$ appearing in \eqref{n1} and \eqref{n2} is called  the {\it Nevanlinna measure}. \\

In a couple of last decades, representations of the form \eqref{n1} or \eqref{n2}, with $\bbr=0$, appeared in the  \emph{free-probability}  as the free-probability analog of the classical  L\'evy-Khintchine formula for infinitely divisible characteristic functions (Fourier transforms).  More precisely, they  are defined as
\begin{equation}\label{voicd}
V(z):=\aar+\int_{\Rset}\frac{1-zx}{z+x}\rho(dx), \quad z\in\Cset\setminus\Rset
\end{equation}
and are called {\it Voiculescu transforms}, cf.  Bercovici and Voiculescu \cite[Section 5, in  particular Theorem 5.10]{bv}. Observe that if $support(\rho)=A\subset \Rset$, then
\begin{itemize}
\item $-V \in \Nev(A)$;
  \item if $\widehat{\rho}$ is the image of $\rho$ by the mapping $x\mapsto -x$, then
  $$z \mapsto V(-z):=\aar+\int_{-A}\frac{zx-1}{z+x}\widehat{\rho}(dx)\in \Nev(-A);$$
  \item  if $\check{A}$ and $\check{\rho}$ are the images of $A\setminus \{0\}$ and $\rho_{A\setminus \{0\}}$ by the mapping $x\mapsto 1/x$, then
$$z \mapsto V\left(\frac{1}{z}\right) =\aar+\rho(\{0\})z +\int_{\check{A}}\frac{zx-1}{z+x}\check{\rho}(dx)\in \Nev(\check{A}).$$
\end{itemize}

Note  that  in \eqref{stars}, we need to know $F$ in some strips  of the complex plane  to retrieve  a measure $\rho$. A natural question is:
\begin{center}\label{natu}
{\it What can be said  about a measure  $\rho$  if we only  have values $F(iw)$, for $w\neq 0,$ and  we\\  don't know if  it is  a restriction of  an analytic function to the  imaginary axis?}
\end{center}
The answer was given by Jankowski and Jurek \cite[Theorem 1]{jan}, where there is  an inversion procedure, which allows to identify the measure $\rho$, or more precisely its  Fourier transform   $\Fou[\rho, s], \; s \in \Rset$, and which is justified as follows: let us, for an index $X$  (where $X$ can be a random variable, or a stochastic process or a measure), define the function  $F_{\Xs}$ on  the  imaginary  axis   $i(\Rset\setminus{\{0\}})$ by
 \begin{equation}\label{(1)}
F_{\Xs}(iw):= \aar_{\Xs} + i \,\bbr_{\Xs}\, w +\int_{\Rset}\frac{iwx-1}{iw+x} \rho_{\Xs}(dx), \quad w\neq 0,
\end{equation}
where $\aar_{\Xs}\in \Rset$, $\bbr_{\Xs}\geq 0$  and  $\rho_{\Xs}$ is a non-negative, finite Borel measure.
Furthermore, if
$$\Fou[\rho_{\Xs}\,;s\,]:= \int_\Rset e^{isx} \rho_{\Xs}(dx), \quad s\in \Rset,$$
is the   Fourier transform  (respectively  the {\it characteristic function}) of a finite  measure (respectively probability measure) $\rho_{\Xs}$, then the \emph{Laplace transform} of $\Fou[\rho_{\Xs}\,;s\,]$ satisfies the equality
\begin{eqnarray}
\mathfrak{L}[\,\Fou[\rho_{\Xs}\,;s\,]\,;\, w]&:=& \int_0^\infty e^{- w s}\, \Fou[\rho_{\Xs}\,;s\,]ds = \int_{\Rset} \frac{1}{w-ix}\, \rho_{\Xs}(dx),\quad w>0\nonumber\\
&=& \left\{ \begin{array}{llcr}
\frac{1}{w^2-1}\, \big [i F_{\Xs}(i w) - i \Re\big(F_{\Xs}(i)\big)+ w\, \Im\big(F_{\Xs}(i)\big)\big] &\mbox{if} &w\neq 1\\
&&\\
\displaystyle \int_{\Rset}\frac{1+ix}{1+x^2}\,\rho_{\Xs}(dx) &\mbox{if} & w=1 .\end{array} \right.
\label{vraie}
\end{eqnarray}

Motivated by the special structure of the Nevanlinna-Pick functions appearing  in \eqref{vraie}, we focus in Section \ref{free} on their analytical aspect and provide additional clarification to their relation to classical and free probability context. An appendix  is included to explain properties of some used special functions and Section \ref{proo} is devoted to the proofs. \\
\medskip

Our main results are obtained in sections \ref{nevaid} and \ref{app}, where we study two  aspects of  Nevanlinna-Pick functions:
\begin{itemize}
  \item The first studied aspect, is mainly probabilistic. In Section \ref{nevaid} we generalize the connection  described by Schilling, Song \& Vondra\v{c}ek \cite[Theorem 6.9]{SSV}, between the subclass of nonnegative functions in $\Nev$ and the one of Bernstein functions associated to subordinators. We provide an  connection of the same nature between a subclass of Nevanlinna-Pick functions and spectrally negative L\'evy processes. This connection is also of same nature than the one obtained right after \eqref{cbf} below when replacing {\it Bernstein functions} by {\it L\'evy-Laplace exponents}. We also provide a property of temporal complete monotonicity for subordinators as Bertoin and Yor \cite{by} did for self-similar Markov processes obtained via the {\it Lamperti transform}, cf. \eqref{temp} below. During several discussions of the first author with Lo\"{i}c Chaumont, a recurrent question was:
\begin{eqnarray} &\mbox{\it  Let $\phi$ be the Bernstein function of some subordinator $\xi$ such that
is impossible to check if $\Psi=\phi^{-1}$}\nonumber\\
&\mbox{\it is a L\'evy-Laplace exponent by standard calculations. Are there  any other clue indicating  whether}\nonumber\\
&\mbox{\it $\xi$ is the first passage time of some spectrally negative L\'evy process?} \label{genui} \end{eqnarray}
      Motivated by this question, we provide in Theorem \ref{temporalcm} a temporal complete monotonicity property:  Any subordinator  $\xi={(\xi_t)}_{t\ge 0}$, such that its $\Psi$-function satisfies some integrability conditions, satisfies one the following:
      \begin{enumerate}[{\it (i)}]
      \item $t \mapsto  \er[\xi_t^{-p}]\;$ is a Bernstein function, for some $\,p \in (-1,0)$;
      \item $t \mapsto  \er[\xi_t^{-p}]\;$ is a completely monotone function, for some $\,p \in [0,1)$;
       \item $t \mapsto t\; \er[\xi_t^{-p}]\;$ is a completely monotone function, for some $\,p \geq 1$.
      \end{enumerate}
      In  Corollary  \ref{tempo}, we answer to question \eqref{genui}  by showing that a subordinator  $\xi$  is the inverse times of a spectrally negative L\'evy process, if, and only if,
      $$t \mapsto t \, \er[\xi_t^{-p}]\;\; \mbox{is a {\it Stieltjes transform}, for some $\,p\geq 1$}.$$
      The latter is also equivalent to:
      $$t \mapsto t^2 \,\er[\xi_t^{-p}]\;\; \mbox{is a {\it complete Bernstein function}, for some $\,p\geq 1$}.$$
      Of course, one can argue that the calculus of the moments $\er[\xi_t^{-p}], \,t>0, \, p>-1,$ are also not always feasible,  but, when they are, we immediately have a conclusion, see Example \ref{examplestable} and especially Example \ref{examplelambert}  . To the best of our knowledge, the temporal properties that we obtain are new.

  \item The second aspect, studied in Section \ref{app}, is related to the free-probability context, and consists in the application  of the inversion formula \eqref{(31)}. In  Propositions \ref{pro1}, \ref{pro2}, \ref{pro4} and Corollary \ref{pro3},  below, we show that for some related index $X$, the functions
      $$w\mapsto \frac{i F_{\Xs}(i w)}{w^2-1},\quad w>1,$$
      in  \eqref{vraie},  are indeed,  Laplace  transforms of some  functions or  measures that, in principle, enables us to identify the corresponding  measure $\rho_{\Xs}$. Let $X=C,\, S$ and $T$  be random variables defined by  their  hyperbolic  characteristic functions $$\Four_{\Xs}(s):=\er[e^{is X}]=\int_{\Rset}e^{is x} \, \pr(X\in dx), \quad  s \in \Rset,$$
$$\Four_{\Cs}(s):=\frac{1}{\cosh(s)},\quad  \Four_{\Ss}(s):=\frac{s}{\sinh(s)}, \quad \Four_{\Ts}(s):= \frac{\tanh(s)}{s}.$$
and  $$\mbox{$\tilde{X}$ be the {\it  free-analog} of the of $X,$ obtained by the procedure {\bf(P4)} in Section \ref{free}}.$$
We emphasize that in our approach to free-probability theory, we only use  purely imaginary numbers $\,z=iw, \,w>1$; cf. \eqref{(31)} and Jurek \cite[Corollaries 3, 4 and 5]{j2} and we recall that the hyperbolic characteristic functions  were studied:
\begin{enumerate}[{\it (i)}]
\item from an infinite divisibility point of view, in Pitman and Yor \cite{Pitman};
\item from a  self-decomposability point of view in  Jurek \cite{j1}, as infinite series of independent exponentially distributed variables;
\item from stochastic representations of their {\it background driving L\'evy processes}  (BDLP) in Jurek-Vervaat \cite[Theorem 3.2]{jv}. The last can be done since all hyperbolic characteristic functions are self-decomposable ones  and therefore, they  admit  a representation  by random integrals
\begin{equation}\label{bdl}
X=\int_0^\infty  e^{-t}\,dY_t,\quad \mbox{where $Y:=Y_{\Xs}$ is the  BDLP},
\end{equation}
see also Jurek and Mason \cite[Chapter 3, Theorem 3.6.8]{jm}.
\end{enumerate}
\end{itemize}
\section{Various interpretations of Nevanlinna-Pick functions}\label{free}
We first list several interesting properties of Nevanlinna-Pick functions.
\medskip

\begin{enumerate}[{\bf(P1)}]
\item Every Nevanlinna-Pick function $F$ in \eqref{n1} is represented by
$$F(z)=\aar + \bbr z+ \int_{\Rset} f(x,z)\,\rho(dx),\quad z\in \Cset\setminus \Rset,$$
where the functions
\begin{equation}
z \mapsto f(x,z)=\frac{xz-1}{z+x}, \quad x\in \Rset,\;z\in \overline{\Cset}\backslash \{-x\}
\label{mob}\end{equation}
are the  {it M\"{o}bius mappings}, which are conformal (i.e analytic with non-null derivative) on any region contained in $\overline{\Cset}\backslash \{-x\}$ and then, preserve the angles in this region. Roughly speaking, $f$ preserve the angles means that for any two rays $L$ and $L'$, starting at a point $z_0$,
the angle which their images $f(L)$ and $f(L')$  make at $z_0$ is the same as that made by $L$ and $L'$, in size as well as in orientation. These mappings enjoy the following properties:
\begin{eqnarray*}
&(a)&  f(x,i) = i; \\
&(b)& \Im \big(f(x,z)\big)= \frac{1+x^2}{|x+z|^2}\Im(z) \; \,\Longrightarrow \;\,  f(x,\Hset^+)\subset\Hset^+ \;\; \mbox{and}\;\; f(x,\Hset^-)\subset  \Hset^-\,; \\
&(c)& f(x,\overline{z})=\overline{f(x,z)};\\
&(d)& f(x,-z)=-f(-x,z);\\
&(e)& f\left(\frac{1}{x},\frac{1}{z}\right)=-f(x,z), \quad x,\,z \neq 0.
\end{eqnarray*}

\item For every $\alpha > 0$ there exists $N_\alpha$  such that $F$ has a left inverse $F^{-1}$
defined on the region  $\Gamma_{\alpha,N_\alpha}$ given by \eqref{reg}. cf. \cite[Proposition 5.4 and Corollary 5.5]{bv}.

\item Any {\em Cauchy transforms}, i.e., a function of the form
  \begin{equation}\label{sti}
   G_\mu(z)= \ccr +  \int_{\Rset} \frac{1}{z+x}\mu(dx), \quad z\in  \Cset {\footnotesize \setminus} \Rset,
   \end{equation}
   where $\ccr$ is a real number and  $\mu(dx)=(1+x^2)\rho(dx)$ is a finite measure, satisfies $-G \in \Nev(\mathds{R})$. To see that, it is enough to take  $\bbr=0$ in \eqref{n2} and to replace $\aar +\int_{\Rset}x \,  \rho(dx)$ by $-\ccr$. In other terms, $G$ swaps the upper and lower complex half planes. In this case, $z\mapsto G_\mu(1/z) \in \Nev(\mathds{R})$. To see this, use property (c) of the function $f$ in \eqref{mob} or the half-plane swapping property, which also yields $H_\mu=1/G_\mu \in \Nev(\mathds{R}) $. Now, assume $\ccr=0$ and observe that
   \begin{equation}\ddr:= \mu(\Rset) = \lim_{|z|\to \infty, \, z\in \Gamma_{\alpha,N}} z \, G_\mu(z)>0,\quad \mbox{for every}\;  \alpha,\;N>0.
   \end{equation}
    The finite measure $\mu$ is called (additive)-{\it free infinitely divisible} if, and only if,  the function $H_\mu=1/G_\mu$ is such that its inverse $H_\mu^{-1}$ (which, by {\bf (P2)} in Section \eqref{free}, always exists) satisfies
   \begin{equation}
   z\mapsto F_\mu(z):=\ddr \,z-H_\mu^{-1}(z) \in \Nev(\mathds{R}).
   \label{freed} \end{equation}
   The function $V_\mu(z)=-F_\mu$ is commonly called the {\it Voiculescu transform} of $\mu$, cf. \eqref{voicd} and Bercovici and  Voiculescu \cite{bv}. The {\it free probability context} corresponds to probability measures $\mu$ and the parameter $\ddr$ in \eqref{sti} is set to be equal 1.

\item Let $\Psi$ be a {\it L\'evy-Laplace  exponent}, i.e., a function represented by
    \begin{equation}\label{EL}
    \Psi(\lambda)= \alpha +  \beta \lambda +  \frac{\gamma^2}{2}  \lambda^2 + \int_{\Rset \setminus \{0\}}\left(e^{-\lambda x}-1+\lambda \frac{x}{1+x^2}\right)\nu(dx), \quad \lambda\in i \,\cdot \, \Rset, \end{equation}
    where $\alpha, \, \gamma \geq 0$, $\beta \in \mathds{R}$,  and the so-called {\it L\'evy measure} $\nu$ is supported by  $\Rset \setminus \{0\}$ and such that
    \begin{equation}\label{nunu}
    \int_{\Rset}\left(x^2 \wedge 1\right) \nu(dx) <\infty \Longleftrightarrow \int_{\Rset} \frac{x^2}{1+x^2}  \nu(dx) <\infty.
    \end{equation}  As we did for the Nevanlinna class, we denote
    \begin{equation}\label{asin}
    \Lle(A):=\big\{\mbox{Laplace exponents}\,\; \Psi  \; \;\mbox{represented by \eqref{EL}, s.t.}\;support(\nu)\subset A  \big\}, \quad  A\subset \Rset.
    \end{equation}
    For the set $\Lle(0,\infty)$, we have the following stochastic interpretation: there is a bijection between the class of  (killed) L\'evy process, i.e. a process $Z={(Z_t)}_{t\geq 0}, \; Z_0=0$,  with stationary and independent increments, and the class of L\'evy-Laplace  exponent $\Psi_{\Xs}$, represented by \eqref{EL},  via the celebrated L\'evy-Khintchine formula, i.e. the Laplace representation: if $X:=Z_1$, then
   $$\er[e^{\lambda Z_t}] = e^{t \Psi_{\Xs}(\lambda)}, \qquad \mbox{for}\;\,t\geq 0, \;\;\lambda \in i.\Rset,$$
    Actually, the distributions of $Z_t,\,t>0$, are entirely determined by the one of
    the {\it infinitely divisible} random variable $X=Z_1$. In  \eqref{EL}, the usage is to call the quantity  $\alpha_{\Xs}$ the {\it killing rate}, $\beta_{\Xs}$ the {\it drift term} and $\gamma_{\Xs}$ the {\it Brownian coefficient}. To every function Nevanlinna-Pick function $F$, represented by \eqref{(1)}, one can associate a free-infinitely divisible random variable $\tilde{X}$ that we call {\it free-analog} of $X$ by these means:
    \begin{itemize}
    \item Pick a Voiculescu transform of some free-infinitely divisible random variable  $\tilde{X}$ whose distribution is  the probability measure $\mu$ in obtained by \eqref{freed}.\\
    \item Barndorff-Nielsen and Thorbj{\o}rnsen \cite[Theorem 4.1]{BT} and  Jurek \cite[Theorem 1]{j2}, clarified, in a more involved way, the so-called Bergovici-Pata bijection from the class of L\'evy processes to the one of free-infinitely divisible random variables, cf.   Bercovici, Pata \& Biane \cite{BP}. More precisely, there exists an infinitely divisible random variable $X$ with L\'evy-Laplace exponent $\Psi_{\Xs}$, represented by \eqref{EL}, such that the function
    \begin{equation}
    F_{\Xts}(iw)= i \, w^2  \int_0^\infty e^{-w u} \,\Psi_{\Xs}(-iu)\,du, \quad w>0. \label{FL}
    \end{equation}
    \item Observe that the Laplace transform of $\Psi_{\Xs}(-iu)$ always exists since it is a continuous function which satisfies  $\lim_{|u|\to \infty}\Psi_{\Xs}(-iu)/u^2=-\gamma^2/2$. The latter is justified by the control of the kernel: for $\Re(\lambda) \geq 0, \; x\neq 0$, we have
        \begin{eqnarray}
        e^{-\lambda x}-1+ \lambda \frac{x}{1+x^2}&=& (e^{-\lambda x}-1+ \lambda x)\,\II_{|x|\leq 1}+  (e^{-\lambda x}-1)\,\II_{|x|\geq 1}   +\lambda \left(\frac{x^2}{1+x^2} \,\II_{|x|\leq 1} - \frac{x}{1+x^2}\,\II_{|x|\geq 1}\right) \label{cont}
        \end{eqnarray}
        and for $\lambda \geq 0$ or $\lambda= iu, \, u \in \Rset$, we have
               \begin{eqnarray*}\label{contp}
        \left|e^{-\lambda x}-1+ \lambda \frac{x}{1+x^2}\right|  &\leq&   \frac{|\lambda|^2 \,x^2}{2}\,\II_{|x|\leq 1}+  2\,\II_{|x|\geq 1} + |\lambda|  \,\left(\frac{x^2}{1+x^2} \,\II_{|x|\leq 1} + \frac{x}{1+x^2}\,\II_{|x|\geq 1}\right) \leq   C(\lambda)\, (x^2 \wedge 1),
        \end{eqnarray*}
         where $C(\lambda):= 3+ |\lambda|+ \frac{|\lambda|^2 }{2}$. The latter insures that  $\Psi_{\Xs}$ is well defined due to the condition \eqref{nunu} on its L\'evy measure. Thus, the triplet of characteristics of $F_{\Xts}$ are given by
    \begin{equation}\label{level}
    \aar_{{\Xts}}= \beta_{\Xs}, \quad \bbr_{{\Xts}}= \alpha_{\Xs} \quad \mbox{and}\quad \rho_{{\Xts}}(dx)= \gamma_{\Xs}^2 \delta_0(dx)+\frac{x^2}{1+x^2}\nu_{\Xs}(dx),
    \end{equation}
    $\rho_{{\Xts}}$ is a finite measure and
    \begin{center}  {\it the operator $\Delta :\Psi \mapsto F$ defined by \eqref{FL}, induces a bijection from $\Lle(0,\infty)$ to $\Nev(\Rset)$, or equivalently, \\from the class of L\'evy processes to the one of free-infinitely divisible random variables.}
    \end{center}
    \begin{remark}\label{simdis} Assume $X:\simdis X_1+X_2$ where $X_1$ and $X_2$ are independent infinitely divisible random variables. On the level of the L\'evy measures, we have $\nu_{\Xs}= \nu_{X_1}+ \nu_{X_2}$ and \eqref{level} gives  $\rho_{\Xts} =\rho_{\Xts_1} +\rho_{\Xts_2}$ on the level of the Nevanlinna measures.
    \end{remark}
    \end{itemize}
\item If $F$ has a continuous extension with $F: (0,\infty)\to \mathds{R}$, then  $F\in \Nev[0,\infty)$. Further, if
    $$F(0+)=\aar -\int_{(0,1]}\frac{1}{x}  \rho(dx)$$ is finite, then $F\in \Nev(0,\infty)$, cf. \cite[Theorem 6.9]{SSV}.

\item  The class $\ST$ of {\it Stieltjes transforms} is formed by functions with a represented similar to \eqref{sti}: $f\in \ST$, if
    \begin{equation}\label{stiel}
     f(\lambda)= \ddr  + \frac{\qqr}{z}+    \,\int_{\oi}  \frac{1 }{z +u} \; \Delta (du), \quad z\in \Cset_+,
    \end{equation}
     where $\ddr, \qqr \geq 0$ are constants and $\int_{\oi}  (1 +u)^{-1} \; \Delta (du) < \infty$, cf.  \cite[Definition 2.1]{SSV}. The class $\ST$ is simply seen as the class of double iterated Laplace transforms, due to the representation
     $$ \frac{1}{1+z} \int_0^{\infty} e^{-z x} \, e^{-x}\,dx, \quad z \geq 0$$
    Using the  the {\it exponential integral function} $\Ei$  defined by \eqref{ei} in the Appendix, an example is given by the function $\ei(z):-e^z \,\Ei(-z)$ represented by
    $$\ei(z):= \int_0^{\infty} \frac{e^{-x}}{z+x}\,dx, \quad z> 0.$$
    The function
    $$\mathfrak{S}(z)=\frac{\log z}{z-1}= \int_{0}^{\infty}\frac{1}{z+1}\, \frac{1}{x+1} \,dx , \quad  z>0,$$
    is another important example of a Stieltjes function and, actually, it is a double iterated Stieltjes function.  As for the class $\ST$, it appears as the kernel of the class s $\ST_2$, of double iterated Stieltjes transforms of functions, whose  several properties were provided by Yakubovich and Martins \cite{YM}. This class could be extended to iterated Stieltjes transforms of measures
    $$\ST_2:=\left\{f(z)= \int_{(0,\infty)} \mathfrak{S}(zt) \;  \Delta(dt), \quad z>0 \right\}.$$
    We believe that it might be interesting to study the properties of the class of iterated Laplace transform of any order, and we hope that this will be the scope of a future work.\\

    By Theorem \cite[Theorem 6.2, Proposition 7.1, Theorem 7.3]{SSV}, the class $\CBF$ of {\it complete Bernstein functions} is formed by those functions $f$ of one the following forms
    $$f(z)=z \,g(z)= h\left(\frac{1}{z}\right)=\frac{1}{k(z)}, \quad \mbox{where}\quad  g, \,h,\, k \in \ST.$$
    Notice that $\CBF$     is contained into the class $\BF$ of {\it Bernstein functions} formed by functions $\phi$ of the form:
    \begin{equation}\label{BF}
     \phi(z)  =\qqr+ \ddr z  + \int_{\oi}(1-e^{-z x})\,\Pi(dx),\quad z\geq 0
    \end{equation}
    where $\qqr\\geq 0$ is the so-called {\it killing rate}, $\ddr \in \mathds{R}$ is the {\it drift} and $\Pi$ is the {\it L\'evy measure} of $\phi$, i.e. a positive measure on $\oi$ which satisfies $\int_{\oi}\left(x \wedge 1\right) \Pi(dx) <\infty$ and that $\phi'(z)$ is a {\it completely monotone function}. We recall that a function $f$ is completely monotone function on $\oi$ if $f$ in infinitely often differentiable and $(-1)^nf^{(n)}\geq 0$, for all $n=0,1,2,\cdotp$. equivalently $f$ is the Laplace transform of some measure, represented by
    \begin{equation}\label{cm1}
    f(\lambda)=\int_{[0, \infty)}e^{-\lambda x} \mu(dx)= \lambda \int_0^\infty e^{-\lambda x} \mu([0,x))\, dx.
    \end{equation}
    Cf. \cite[Chapter 1]{SSV}. Complete Bernstein function
    were defined in \cite[p. 69]{SSV}  as those Bernstein function such that their L\'evy measure $\Pi$
    have the form
    \begin{equation} \label{picb}
    \Pi(dx)= m(x)\, dx,\;>0, \quad \mbox{where $m$ is a completely monotone function},
    \end{equation}
i.e. $m(x)=\mathfrak{L}[\sigma\,;\,x]$. Thus, every $f\in \CBF$ takes the $\Nev(0,\infty)$ form:
    \begin{equation}
    f(\lambda)=  \qqr+ \ddr \lambda+ \int_{(0,\infty)} \frac{\lambda}{\lambda+x}  \sigma(dx), \quad \lambda \geq 0,\label{cbf1}
    \end{equation}
    where $\qqr,\,\ddr \geq 0$ and $\sigma$ is some measure on $(0,\infty)$ that integrates $1/(1+x)$.
    From the proof of \cite[Theorem 6.9]{SSV}, one can extract  that the Nevanlinna-characteristics of $f$ are
    \begin{equation}\label{chn}
    \aar_f=  \qqr+ \int_{(0,\infty)}\frac{\sigma(dx) }{1+x^2}, \quad \bbr_f= \ddr \quad \mbox{and}\quad \rho_f(dx)=\frac{x}{1+x^2}\sigma(dx).
    \end{equation}
    By \cite[Theorem 6.2]{SSV}, the extension of $f$ is also represented by a Nevanlinna representation analogous of \eqref{FL}:
    \begin{equation}
    f(z)=z^2  \int_0^\infty  e^{-z x} \,\phi(x)\,dx= z^2  \, \Lap[\phi,z], \quad z \in \Cset^+, \label{cbf}
    \end{equation}
    where $\phi$ is a some Bernstein function. The latter induces that
    \begin{center}  {\it the operator $\phi \mapsto f$ given by \eqref{cbf}, is a bijection from $\BF$ to $\CBF$, or equivalently, from the class \\of subordinators to the one of complete subordinators.}
    \end{center}
    By \cite[Theorem 6.9]{SSV}, we extract that actually, the analytic extension on $\Cset\setminus \Rset_+$ of $\CBF$-functions  corresponds to those functions  $F$ in $\Nev(0,\infty)$ which are nonnegative-valued on $(0, \infty)$. A non trivial example of a $\CBF$-function is given by the following: let $\Li_s$  be the  {\it polylogarithm functions} given by \eqref{pl} and define $$\li_s(z):= -\Li_s(-z), \quad z\geq 0, \; \, \;s>0.$$
    By \eqref{(81)}, $\li_s$ meets the $\CBF$-form \eqref{cbf1}:
    $$\li_s(z)=\frac{1}{\Gamma(s)}\int_1^\infty  \frac{z}{z+x} \,\log^{s-1}(x)\, \frac{dx}{x}, \quad z\geq 0.$$
    If furthermore $s\geq 1$, $\li_s$ is in the class $\TBF$ of {\it Thorin Bernstein functions}, the subclass of $\CBF$-functions, such that the $m$-function in \eqref{picb}
    \begin{equation} \label{pitb}
    \Pi(dx)= m(x)\, dx,\;x>0, \quad \mbox{where $x\,m(x)$ is a completely monotone function},
    \end{equation}
    cf. \cite[Theorem 8.2 (v)]{SSV}. In general, $\li_s, \, s>0,$   meets the Nevanlinna $\Nev[1,\infty)$-form \eqref{n1} on $\Cset {\footnotesize \setminus} [-1,\infty)$ and its characteristics are given by \eqref{chn}:
    $$\aar=\frac{1}{\Gamma(s)}\int_1^\infty  \frac{\log^{s-1}(x)}{x(1+x^2)}\,  dx, \quad  \bbr=0  \quad \mbox{and}\quad \rho(dx)=\frac{1}{\Gamma(s)}\,\frac{\log^{s-1}(x)}{1+x^2}\, dx, \; \;x\ge 1.$$
    By formula \cite[(1.110)  p. 27]{Lewin}, observe that $\li_2$ and $\ei$ are linked by a formula of type \eqref{cbf}:
    $$z\mapsto z^2  \int_0^\infty e^{-z x} \,\li_2(x) \,dx= z \int_z^\infty \frac{\ei(y)}{y}\, dy \, \in \CBF.$$
\item The class $\CBe$ of {\it extended complete Bernstein functions} is formed by those functions $f$ of the form $f=g-h$, where $g\in \CBF$ and $h\in \ST$.   Actually, the analytic extension on $\Hset^+$ of $\CBe$-functions corresponds to those   functions  $F$ in $\Nev[0,\infty)$ which are real-valued on $(0, \infty)$. See \cite[Proposition 6.12]{SSV} and Remark \ref{coco} below.
\end{enumerate}
After these preliminaries, we are able to state our first result which is a generalisation of the inversion procedure \eqref{vraie} and of \cite[Corollary 6.13]{SSV} stated for $\CBe$-functions.
\begin{theorem} Let $A\subset \Rset$ and $F:\Hset^+\to \Hset^+$ be a an analytic function and recall $\mathfrak{L}$ and $\Fou$ stand  of the Laplace and Fourier transforms. Then, the following statements are equivalent:
\begin{enumerate}[1)]
\item $F\in \Nev(A)$;
\item $\displaystyle G(z):=\frac{1}{1+z^2} \, \left[\Re\big(F(i)\big) + \Im\big(F(i)\big)\,z- F(z)\right], \; z\in \Hset^+, \;$  is the Cauchy transform of a finite measure on $A$.
\item The function $\displaystyle H(w) :=\frac{iF(iw)}{w^2-1}, \;\, w>1$, is of the form
\begin{equation}
H(w)= \mathfrak{L}\left[\Fou[\rho\,;\,s]+i \Re\big(F(i)\big) \sinh(s) - \Im\big(F(i)\big) \cosh(s); w\right].
\label{(31)}\end{equation}
for some finite measure $\rho$ on $A$.
\end{enumerate}
In particular $F\in \Nev(\Rset_+),\;$  if, and only if the function $G$ has continuous extension $\;G: (0,\infty) \to \Rset_+$ which is the Stieltjes transform of some finite measure on $\Rset_+$.
\label{vraiz}\end{theorem}
\section{The $\Theta$-transform, Nevanlinna-Pick functions and L\'evy processes} \label{nevaid}
A (possibly killed) subordinator $\xi={(\xi_t)}_{t\geq 0}$ started from 0 is an increasing L\'evy process whose distribution are also obtained  via a L\'evy-Khintchine formula
\begin{equation}
\er[e^{-\lambda \xi_t}]= \int_{[0,\infty)}e^{- \lambda x} \, \pr(\xi_t\in dx)= e^{-t \,\phi(\lambda)}, \quad \mbox{ for} \;t,\;\lambda \geq 0.
\end{equation}
where $\phi$ is  Bernstein function, i.e. a function represented by \eqref{BF}. See the monograph of Bertoin \cite{bertoin}  for background on subordinators and the book of Schilling, Song \& Vondra\v{c}ek \cite{SSV}  \cite{SSV} for Bernstein functions. Lamperti's transformation uses the implicit time-change $\tau_t, \; t\geq 0$, defined by the identity
$$t=\int_0^{\tau_t} e^{\xi_s}\,ds,$$
where $\xi$ is a general subordinator. The process $E={(E_t)}_{t\geq 0}:= {(e^{\xi_{\tau_t}})}_{t\geq 0}$ is a strong Markov process started from 1  which enjoys the scaling property: if $\mathbb{P}_x,\; x > 0$, is the distribution of the process $x\,E_{t/x},\; t \geq 0$,
then $\mathbb{P}_x$ coincides with the law of the process $E$ started from $x$ (that is, when the subordinator $\xi$ is replaced by $\xi + \log x$ in Lamperti's transform). In \cite{by}, it was shown that for all $ p>0$, the function $t\mapsto \er[E_t^{-p}], \;t\geq 0$, is completely monotone  with Laplace transform representation
\begin{equation}
\er[E_t^{-p}]= \int_{[0, \infty)} e^{-t x} \,\sigma_p(dx),
\label{temp}\end{equation}
where the entire moments of the probability measure $\sigma_p$ are given by
$$\int_{[0, \infty)} x^k \,\sigma_p(dx)= \phi(p)\,\phi(p+1)\cdots \phi(p+k-1), \quad k=1,2\cdots$$

A proper subclass of subordinators is formed by first passage times of {\it spectrally negative L\'evy processes} $Z={(Z_t)}_{t\geq 0}$. That means that $Z$ is a L\'evy process with no positive jumps and the  distribution of $Z_t$ is obtained by a L\'evy-Laplace exponent $\Psi$ such that $\er[e^{\lambda Z_t}] =e^{t\Psi(\lambda)}$ and in \eqref{EL}, the L\'evy measure $\nu$ of $\Psi$ is supported by $\oi$. This entails that $\Psi$ has the continuous extension on $\Rset_+$:
\begin{equation}\label{LE}
\Psi(\lambda)= \alpha +\beta \lambda+  \frac{\gamma^2}{2} \lambda^2 + \int_{\oi}\left(e^{-\lambda x}-1+\lambda \frac{x}{1+x^2}\right)\nu(dx),\quad \lambda \geq 0. \end{equation}
i.e. $\Psi \in \Lle(0,\infty)$. Observe that if  in \eqref{LE},     $\nu$ satisfies
\begin{equation}\label{inte}
\overline{\beta} :=  \beta +\int_{\oi} \frac{x}{1+x^2} \nu(dx) \quad \mbox{is a finite quantity,}
\end{equation}
 then $\Psi$ satisfies
\begin{equation}\label{formf}
\lambda \mapsto \varphi(\lambda):= \alpha + \overline{\beta} \lambda + \frac{\gamma^2}{2}\lambda^2 -\Psi(\lambda)  =  \int_{\oi} (1-e^{-\lambda x})\nu(dx) \in \BF.
\end{equation}
Conversely, every Bernstein function $\phi$ satisfies
$$\lambda \mapsto \alpha +  \beta \lambda + \frac{\gamma^2}{2}\lambda^2-\phi(\lambda) \in \Lle(0,\infty), \qquad \mbox{for all}\; \alpha, \, \gamma \geq 0, \; \beta \in \Rset.$$

By the proof of \cite[Theorem 3.12 and Section 8.1]{kyp}, under the
assumption that
\begin{equation}\label{assum}
Z\;\mbox{is non-killed, (i.e.} \;\alpha=\Psi(0) =0) \quad\mbox{and} \quad \er[Z_1]=\Psi'(0+)=\beta -  \int_{\oi}  \frac{x^3}{1+x^2}\,\nu(dx) \geq 0,
\end{equation}
we have unique (nonnegative) solution $\phi(\mu)$   to the equation $\Psi(\lambda) = \mu$, given by the right inverse
\begin{equation}\label{RI}
\phi(\mu) = \sup\{\lambda \geq 0 \,: \, \Psi (\lambda) = \mu\}.
\end{equation}
Necessarily $\phi(0)=0$, and $\phi$ is a Bernstein function  represented as in \eqref{BF}. Actually $\phi$ corresponds to a subordinator $\xi$  which is the first passage of $Z$ time above a level $s > 0$, i.e.  \
\begin{equation}\label{rig}
\xi_s  = \inf\{t > 0 \,:\, Z_t > s\}.
\end{equation}

\medskip
Due to \eqref{cont}, observe that for every $\Psi \in  \Lle(\Rset_+)$,  there exists $A, \, B, \, C \geq 0$, such that
$$\left|\Psi(\lambda) \right|\leq  A + B \lambda + C \lambda ^2, \quad \mbox{for all $\lambda \geq 0$},$$
and then the Laplace transform of $\Psi$ is well defined. With this observation together with \eqref{formf} and even if the integrability condition \eqref{inte} is not satisfied, we still be able to establish a correspondence between the Nevanlinna  class $\Nev(\Rset_+)$  with the class $\Lle(\Rset_+)$ given by \eqref{asin} as obtained the one obtained in \eqref{cbf}:
\begin{proposition}[The $\Theta$-transform] To every pair $(\bbr, \Psi) \in \Rset_+\times \Lle(0,\infty)$, is associated a (unique) Nevanlinna $\Nev(\Rset_+)$-function obtained by the operator $\Theta$ given, for all $z\in \Cset^+$, by
\begin{equation}\label{fz}
\Theta(\bbr, \Psi)(z)= F(z):= z^2\int_0^\infty e^{-z u}\,\big (\bbr+\Psi(0)-\Psi(u)\big) \,du=\aar +\bbr z+\int_{[0,\infty)} \frac{z x -1}{z+x}\rho(dx),
\end{equation}
and the characteristics $(\aar, \rho)$ of $F$ are related to  the one's of $\,\Psi,\,$ $(\alpha, \, \beta, \, \gamma,\, \nu)$, by
$$\aar = -\beta,  \quad   \mbox{and}\quad \rho(dx)= \gamma^2 \,\delta_0(dx) + \frac{x^2}{1+x^2} \nu(dx).$$
Furthermore, the Nevanlinna  class of functions
$$\Nev^\star(\Rset_+):=\left\{F \in \Nev(\Rset_+) \;\;\mbox{s.t.}\; \;  \bbr=\lim_{|z|\to \infty,  \; z\in \Cset^+} \frac{F(z)}{z} =0\right\},$$
is obtained as the image  of $\{0\}\times \Lle(0,\infty)$ by the operator $\Theta$.
\label{zfz}\end{proposition}
\medskip

\begin{remark}\label{coco} As in {\bf (P7)} in Section \ref{free}, one can see that the class
$$\CBe^\star:= \{f \in \CBe, \; \mbox{such that}\;\lim_{\lambda\to \infty} f(\lambda)/\lambda=0\}$$
is such that the  extension on $\Hset^+$ of all  its members $f$  coincides with the class $\Nev^\star(\Rset_+)$ given in last proposition.
\end{remark}
\medskip
\begin{example}[The $\Theta$-transform for spectrally negative stable processes] Let $\vartheta\in (1, 2]$. The function $\Psi_\vartheta(\lambda):=\lambda^\vartheta, \;\lambda\geq 0$, multiplied by positive constants, corresponds to the L\'evy-Laplace exponents of spectrally negative stable processes, see \cite[Section 6.5.3]{kyp}.  Observe that $\Psi_2$ corresponds to the Gaussian distribution and that $\Theta(0, y\mapsto u^2)(z)=-2/z^2$.
For $\vartheta\in (1, 2)$, let
$$c_\vartheta=\frac{\vartheta(\vartheta-1)}{\Gamma(2-\vartheta)}  \quad \mbox{and}\quad \rho_\vartheta(x)=\frac{c_\vartheta\, x^{1-\vartheta}}{1+x^2},\;x>0. $$
Using  \cite[Formula FI II 718 p. 319]{grad} and Euler's reflection formula, we easily get that
$$\beta_\vartheta\:=c_\vartheta\int_0^\infty \frac{x^{2-\vartheta} }{1+x^2}\,dx = \frac{\vartheta}{\Gamma(2-\vartheta)} \Gamma\left(\frac{\vartheta+1}{2}\right) \Gamma\left(\frac{3-\vartheta}{2}\right).$$
The function $\Psi_\vartheta$ is represented by
$$\Psi_\vartheta(\lambda)=  \int_0^\infty  \left(e^{-\lambda x}-1+\lambda x\right) \frac{c_\vartheta}{x^{\vartheta+1}}\, dx =\beta_\vartheta \lambda+  \int_0^\infty  \left(e^{-\lambda x}-1+\lambda\frac{x}{1+x^2}\right) \frac{c_\vartheta}{x^{\vartheta+1}} \,dx.$$
The Nevanlinna-Pick function $F_\vartheta(z)=- \Gamma(\vartheta +1)\,z^{1-\vartheta}, \; z\in \Cset^+$, is the $\CBe$-function   obtained by the operator $\Theta$ in \eqref{fz}. Actually $-F_\vartheta$ is a Stieltjes function in the sense of {\bf (P6)} in Section \ref{free}  and without effort, we get the Nevanlinna representation
$$  F_\vartheta(z)= -\beta_\vartheta + \int_0^\infty \frac{zx-1}{z+x} \, \rho_\vartheta(x)\,dx, \quad z\in \Cset^+. $$
\label{thest}\end{example}
\medskip

Lemma \ref{keylemma} in Subsection \ref{useful} below, is the key result for the following temporal completely monotone property for subordinators.
\begin{proposition}[Temporal completely monotonicity property for subordinators] \label{prop-moment}  Let $\phi:\Rset_+\to \Rset_+$  be a strictly  increasing function such that $\phi(0)=0$ and $\phi(\infty)=\infty$. Let $\Psi$   be its inverse function and $\,\Ebf\,$ denotes a r.v. with standard exponential distribution.  Then  the following statements are equivalent.
\begin{enumerate}[1)]
\item $\phi$ is a Bernstein function;

\item There exists a   subordinator $\xi={(\xi_t)}_{t\geq 0}$  satisfying the identities in law
\begin{equation} \frac{\Ebf}{\xi_t} \simdis  \Psi \left( \frac{\Ebf}{t} \right),\quad \forall t>0,\qquad \mbox{ where $\Ebf$  in the l.h.s.  is assumed to be independent of $\xi_t$};
\label{id2} \end{equation}

\item There exists a family of positive r.v.'s $\eta={(\eta_t)}_{t\geq 0}$  such that  \eqref{id2} holds true with $\xi$ replaced by $\eta$;

\item  There exists a  subordinator $\xi={(\xi_t)}_{t\geq 0}$, such that  the representation
\begin{equation}
 \er[\xi_t^{-p}] = \frac{t}{\Gamma(p+1)}\int_0^\infty e^{-t x} \, \Psi^p(x) \,dx
\label{formulei} \end{equation}
holds true for every $t>0$ and $p>-1$, such that $\Lap[\Psi^p\,;\,t]<\infty$.

\end{enumerate}
Under  the latter, $\phi$ is necessarily the Bernstein function of $\xi$.
\end{proposition}
\begin{remark}
Observe that the temporal completely monotone property in \eqref{formulei} involves the Laplace representation in the right hand side term  similar to the one obtained in \eqref{temp} with $p=1$. Also observe that $\phi(\lambda)=\lambda$ is associated to the trivial subordinator $\xi_t=t$ , then replacing by $\Psi(\lambda)=\lambda$ in \eqref {formulei}, we retrieve the integral representation of the Gamma function \eqref{gag}.
\end{remark}

\begin{example}[Moments of positive stable distributions]
If $\xi^{(\as)}={(\xi^{(\as)}_t)}_{t\geq 0}$ is a standard $\alpha$-stable subordinator, i.e., with Bernstein function $\phi_\alpha = \lambda^\alpha$,  $0<\alpha <1$, then \eqref{id2} gives the following well known results, valid for all $t>0, \;p>-\alpha$: by \eqref{id2} , we have
\begin{equation}\label{cy}
\frac{\Ebf}{\xi^{(\as)}_t} \simdis \left( \frac{\Ebf}{t} \right)^{\frac{1}{\alpha}} \Longrightarrow
\frac{\Ebf}{\xi^{(\as)}_1} \simdis   \Ebf^{\frac{1}{\alpha}} \quad \mbox{and}\quad \xi^{(\as)}_t\simdis t^{\frac{1}{\alpha}}\,\xi^{(\as)}_1\,\;\mbox{({\it scaling property for stable processes}),}
\end{equation}
and by \eqref{formulei} , we have
$$ \er[(\xi^{(\as)}_t)^{-p}] =  \frac{t}{\Gamma(p+1)}\int_0^\infty e^{-t x} \, x^{\frac{p}{\alpha}}\,dx = \frac{t}{\Gamma(p+1)} \, \frac{\Gamma(1+\frac{p}{\alpha})}{t^{1+\frac{p}{\alpha}}}= t^{-\frac{p}{\alpha}}\, \frac{\Gamma(1+\frac{p}{\alpha})}{\Gamma(1+p)},$$
a formula that  in  \cite[Exercise 4.17]{CY}, whereas the second identity in  \eqref{cy} is in \cite[Exercise 4.19]{CY}.
See \cite[Theorem 1]{jedidi} for all moments of all stable distributions.
\label{examplestable}
\end{example}
\begin{remark}\label{simple} Proposition \ref{prop-moment} shows that if $p>-1$ and $\Psi$ is the inverse  of a Bernstein function $\phi$ associated to the subordinator $\xi$, then the following assertions are equivalent.
\begin{enumerate}[(i)]
  \item $t\mapsto \er[\xi_t^{-p}]/t$ is completely monotone on $(0,\infty)$
  \item $\Lap[\Psi^p\,;\,]$ is well defined on $(0,\infty)$.
 \end{enumerate}
\end{remark}
This simple remark is improved as follows:
\begin{theorem}[Temporal completely monotonicity property for subordinators improved] \label{temporalcm}  Let $p>-1$ and  $\phi$  be a non-trivial Bernstein function (i.e. $\phi$ is not affine) such that $\phi(0)=0,\;\phi(\infty)=\infty$. Let $\Psi$ and $\xi={(\xi_t)}_{t\geq 0}$ be respectively, the inverse function and  the  subordinator, associated to $\phi$.  We have the following results.
\begin{enumerate}[1)]
\item If $p\in(-1,0)$, then $e^{-tx}\,\Psi^p(x)$ is integrable near $\infty$ and the following assertions are equivalent.
     \begin{enumerate}[(i)]
     \item $\Psi^p$ is integrable at the neighborhood of $\,0$;
    \item $t\mapsto \,\er[\xi_t^{-p}]$ is a Bernstein function.
    \end{enumerate}
    We then have the representation
    \begin{equation}
    \er[\xi_t^{-p}] = \frac{1}{\Gamma(p+1)}\int_0^\infty (1-e^{-t x}) \,(-\Psi^p)' (x) \,dx
    \label{formula1} \end{equation}
\item If $p\geq 0$, then $\Psi^p$ is integrable near $\,0$ and the following assertions are equivalent.
     \begin{enumerate}[(i)]
     \item $x\mapsto  e^{-tx}\,\Psi^p(x)$ is integrable near $\infty$ for all $t>0$;
     \item $t\mapsto  \er[\xi_t^{-p}]$  is completely monotone on $(0,\infty)$.
     \end{enumerate}
     We then have the representation
    \begin{equation}\label{cmety}
    \er[\xi_t^{-p}] =\frac{1}{\Gamma(p+1)} \int_0^\infty e^{-tx}\,(\Psi^p)' (x) \, dx,\quad t>0.
    \end{equation}
     The function $t\mapsto t\,\er[\xi_t^{-p}]$ is Bernstein if, and only if, $\lambda \mapsto \lambda^{1-p} \phi'(\lambda)$ is non-decreasing, hence  $p\in [0,1)$. In this case, we have the representation
     \begin{equation}\label{formula2}
    t \er[\xi_t^{-p}] =\frac{1}{\Gamma(p+1)} \int_0^\infty (1-e^{-tx})\; (\Psi^p)'' (x) \, dx,\quad t>0.
    \end{equation}
\item If $p\geq 1$, then  $\Psi^p$ is convex,  $(\Psi^p)'(0+)\in [0,\infty)$, and the following two assertions are equivalent.
      \begin{enumerate}[(i)]
      \item $x\mapsto \,e^{-tx}\,\Psi^p(x)\;$ is integrable near $\infty$ for all $t>0$;
      \item $t\mapsto \,t\, \er[\xi_t^{-p}]\;$  is completely monotone on $(0,\infty)$.
     \end{enumerate}
      In this case, we have the representation
     \begin{eqnarray}
      t\, \er[\xi_t^{-p}] &=& \frac{1}{\Gamma(p+1)}\left[(\Psi^p)'(0+) + \int_0^\infty e^{-tx}\, (\Psi^p)''(x) \, dx \right],\quad t>0,\qquad \label{rerp1}
     \end{eqnarray}
     and the following holds true:
      \begin{enumerate}[(a)]
     \item  The function $t\mapsto t^2\, \er[\xi_t^{-p}]$ is Bernstein  if, and only if, $(\Psi^p)'$ is concave;
     \item If further  $\Psi^p\in \Lle(0,\infty)$ with quadruple of characteristics $(0,\,\beta_p, \, \gamma_p, \,\nu_p)$ in its representation \eqref{LE}, then
     \begin{equation}\label{phio}
      (\Psi^p)'(0+)=\beta_p -  \int_{\oi}  \frac{x^3}{1+x^2}\,\nu_p(dx) \in [0,\infty),
     \end{equation}
     and the function $t\mapsto \,t\, \er[\xi_t^{-p}]\;$ has the Stieltjes transform form \eqref{stiel}, i.e.
     \begin{equation}\label{repp4}
      t\, \er[\xi_t^{-p}] = \frac{1}{\Gamma(p+1)}\left[ (\Psi^p)'(0+)+\frac{\gamma_p^2}{t}+\int_{(0,\infty)} \frac{u^2}{t+u}\,\nu_p(du)\right],\quad t>0.
     \end{equation}
     \end{enumerate}
\end{enumerate}
\end{theorem}
\begin{remark}   Under the assumptions of point 3) in Theorem  \ref{temporalcm}, we have
     $$\lim_{t\to \infty} t\, \er[\xi_t^{-p}] =  \frac{(\Psi^p)'(0+)}{\Gamma(p+1)},$$
     and under the assumptions of point 4) in Theorem  \ref{temporalcm}, we have
     $$\lim_{t\to \infty} t^2\, \er[\xi_t^{-p}] -  t\frac{(\Psi^p)'(0+)}{\Gamma(p+1)} = (\Psi^p)''(0+)=\gamma_p^2 +      \int_{(0,\infty)} u^2\, \nu_p(du).$$
\end{remark}
\begin{example}[Moments of Lambert distributions] Besides the example of stable subordinators given in Example \ref{examplestable}, we propose the following less trivial  one:
in Pakes \cite{pakes}, the  principal Lambert function $\phi_{\Ls}$,  defined as the unique real-valued concave increasing solution to the functional equation $\lambda=\phi_{\Ls}(\lambda)\, e^{\phi_{\Ls}(\lambda)},\; \, \lambda\geq 0,$ is shown to be a  Thorin Bernstein, i.e., a Bernstein function with L\'evy measure as in \eqref{pitb}. Let $\xi^{\Ls}={(\xi^{\Ls}_t)}_{t\geq 0}$ be its associated subordinator. Since $\Psi_{\Ls}(x)=x\, e^{x}, \; x\geq 0$, then r.v.'s $\xi^{\Ls}_t$ satisfy the identity \eqref{id2}:
$$\frac{\Ebf}{\xi^{\Ls}_t} \;\simdis \; \frac{\Ebf}{t}\; e^{\frac{\Ebf}{t}}\,  \simdis \; -\frac{1}{t}\; \Ubf^{-\frac{1}{t}}\; \log \Ubf, \quad \mbox{for all $t>0$},$$
where $\Ubf$ has the uniform distribution on $(0,1)$, and also \eqref{formulei}:
\begin{equation} \label{formulei1}
 \er[(\xi^{\Ls}_t)^{-p}] = \frac{t}{\Gamma(p+1)}\int_0^\infty e^{-t x} \; (x \; e^{x})^p \,dx=
\frac{t}{(t-p)^{p+1}},\qquad \mbox{for all}\; t>0,\;\; t>p>-1.
\end{equation}
Thus, for all $t>0$, the function
$$ q\mapsto M_t(q):= t \, (t+q)^{q-1},\qquad - t<q <1,$$
is a Mellin transform whose domain of definition  can be extended to the interval $(-q,\infty)$, i.e., the representation
$$M_t(q)= \er[(\xi^{\Ls}_t)^{q}]=t \, (t+q)^{q-1}\; \, \mbox{remains valid if $q>-t$}.$$
For more arguments  justifying this extension on the Mellin transform, we refer to  Jedidi et al. \cite{JBH} for instance. Note that  Pakes \cite[Theorem 3.4]{pakes},  calculated only the moments of natural numbers order for $\xi^{\Ls}_1$:
$$\er[(\xi^{\Ls}_1)^{q}]=  (1+q)^{q-1},\qquad \mbox{for}\; \,q=0,1,2,\cdots $$
We can go beyond by observing that condition 3)(i) in Theorem \ref{temporalcm}   fails for $\Psi_{\Ls}$. Nevertheless, since $\Lap[(\Psi_{\Ls})^p, p+t]$ is finite for every $p,\,t>0$, then \eqref{formulei1} can be restated
$$t\mapsto (t+p)\,\er[\xi_{t+p}^{-p}] = \frac{t+p}{t^{p+1}}=\frac{1}{t^{p}}+ \frac{p}{t^{p+1}}\; \; \mbox{is a completely monotone function.} $$
\label{examplelambert}\end{example}
\medskip

Let $\phi$ be a Bernstein function, with inverse $\Psi$ and associated subordinator $\xi$. Observe that it is not always possible to have $\Psi$  explicitly, for instance, take the Bernstein function $\phi(\lambda)=\lambda/\log(1+\lambda), \; \lambda \geq 0$. Even one is lucky enough to have
$\Psi$ explicitly, it might be impossible to easily check whether it is a L\'evy-Laplace exponent! For instance, no standard calculation would exhibit that  $\Psi(\lambda)=\lambda/[(1+\lambda)\log(1+1/\lambda)], \; \lambda \geq 0$, is genuinely a L\'evy-Laplace exponent, cf. \cite[Example 40 p. 320]{SSV}.

An answer to  question \eqref{genui} might be a consequence of Theorem \ref{temporalcm}: we show in next result that we only need to strengthen the equivalence in Remark \ref{simple},  $\xi$  needs to have  a {\it temporal Stieltjes property}! Representation \eqref{foms} below is actually connected to the Laplace and Nevanlinna representations in \eqref{fz}. For
\medskip

\begin{corollary}[On temporally Stieltjes property for spectrally negative L\'evy processes] $\;$ \label{corollary1}
\begin{enumerate}[1)]
\item Let $\Psi$ be a L\'evy-Laplace exponent represented by \eqref{LE} and  satisfying
\eqref{phio} wiz.  $\psi'(0+)\in [0,\infty)$.
Let the inverse function of (the non-negative) function $\Psi$ be the Bernstein function $\phi$ associated to the subordinator $\xi={(\xi_t)}_{t\geq 0}$ given by
\eqref{RI} and \eqref{rig} respectively. Then, the following  holds:
\begin{enumerate}[(i)]
\item For all $t>0$, $\xi_t,$ satisfies the identity in law \eqref{id2}.
\item If $\,p\geq 0,\,$  or if either $\,0 > p>-1\,$ and $\,\Psi^p\,$ is integrable at 0, then the function $t\, \,\er[\xi_t^{-p}], \;t>0$, is well defined, is completely monotone, and is represented by \eqref{formulei}.
\item If $\,p\geq 1$ and  $\Psi^p \in \Lle(0,\infty)$,  with quadruple of characteristics  $(0, \, \beta_p, \,\gamma_p, \,\nu_p)$ in its representation \eqref{LE} and satisfying \eqref{phio}, then
      \begin{equation}\label{conclus0}
    t \mapsto   \frac{t^2}{\Gamma(p+1)}\,\Lap[\Psi, t]  =   t  \;\er[\xi_t^{-p}]\;\; \mbox{is a Stieltjes transform,}
       \end{equation}
    and we have the representation
    \begin{equation}\label{foms}
   \Gamma(p+1)\; t  \;\er[\xi_t^{-p}]=\Psi'(0+) +\frac{\gamma^2}{t} +\int_{\oi}
    \frac{x^2}{t+x}\,\nu_p(dx).
    \end{equation}
     The latter is equivalent  to
    \begin{equation}\label{conclus}
    t \mapsto \frac{ t^3}{\Gamma(p+1)} \,  \Lap[\Psi, t] = t^2 \;\er[\xi_t^{-p}]  \quad \mbox{is a complete Bernstein function}.
    \end{equation}
\end{enumerate}
\item  The  converse is stated as follows. Let $p\geq 1$ and ${(\xi^{(\scriptstyle 1/p)}_t)}_{t\geq 0}$ be a $1/p$-standard stable subordinator, (see
Example \ref{examplestable}, with the convention $\xi^{(\scriptstyle 1)}_t=t$). If a subordinator $\xi$ is such that the function $t\mapsto t \;\er[\xi_t^{-p}]$ is a Stieltjes transform,
then the subordinated process ${\big(\xi^{(\scriptstyle 1/p)}\circ  \xi)_t\big)}_{t\geq 0}$ is   the first passage time process of a spectrally negative L\'evy processes with L\'evy-Laplace exponent  $\Psi^p \in \Lle(0,\infty)$.
\end{enumerate}
\label{tempo}\end{corollary}
\begin{remark}  \label{remel1} Last results merit  the following comments:
\begin{enumerate}[{\it (i)}]
   \item The equivalence between \eqref{foms} and \eqref{conclus} is immediate by property {\bf (P6)} above.
  \item Recall the Lambert Bernstein function $\phi_{\Ls}$ given in Example \ref{examplelambert} and its associated subordinator $\xi^{\Ls}$. Since the inverse function  $\Psi_{\Ls}(x)=x\, e^{x}$ is not a L\'evy-Laplace exponent, then $\xi^{\Ls}$ is certainly not the inverse time process of a L\'evy process. Thus, Corollary  \ref{tempo} does not apply.
  \item Assume we know that a subordinator $\xi={(\xi_t)}_{t\geq 0}$ is such  that $\er[\xi_t^{-p}]$ is finite and is explicit for all $t>0$ and some $p>0$, for example $\xi=\xi^{(\scriptstyle \alpha)}$ where is an $\alpha$-stable subordinator. A ``good test'' to check whether $\xi$  is the inverse time of some spectrally negative L\'evy process with L\'evy exponent $\Psi$, is ``{\it just check if the temporal property \eqref{conclus0} or \eqref{conclus} of $\xi$ holds true. If it fails, then $\xi$ is certainly not  an inverse time}''.
  \item Additionally, the operator applied on $\Psi$ in \eqref{conclus} has to be compared with the one obtained by Schilling, Song and  Vondra\v{c}ek \cite[Theorem 6.2.]{SSV} and illustrated by \eqref{cbf}. Recall for instance  the function $\Psi_\vartheta(\lambda):=\lambda^\vartheta,  \; \lambda \geq 0$  used in Example \eqref{thest}.  $\Psi_\vartheta$ is a L\'evy-Laplace exponent if, and only if, $\vartheta\in[1,2]$ and then, $\Psi_\vartheta$  corresponds  to a (non-trivial) spectrally negative stable processes $Z^\vartheta={(Z^\vartheta_t)}_{t\geq 0}$ if, and only if $\vartheta\in(1,2]$. Denoting $\alpha=1/\vartheta$, the  subordinator $\xi^{(\as)}={(\xi^{(\as)}_t)}_{t\geq 0}$  is the inverse in time of $Z^\vartheta$ and is associated to the Bernstein function $\lambda \mapsto \lambda^\alpha, \, \alpha \in [1/2, 1)$.  Observe that $\xi^{(\as)}$ has the so-called temporal {\it scaling property}: $\xi^{(\as)}_t\simdis t^{1/\alpha}\, \xi^{(\as)}_1$. Example \ref{examplestable}  gave, that for $p>-\alpha$,
    $$\er\left[(\xi^{(\as)}_t)^{-p}\right]= c\,t^{-\frac{p}{\alpha}}, \quad c:= \frac{\Gamma(\frac{p}{\alpha}+1)}{\Gamma(p+1)}$$
    and then
    $$t \mapsto t^2\, \er\left[(\xi^{(\as)}_t)^{-p}\right]= c \, t^{2-\frac{p}{\alpha}}\in \CBF\quad \Longleftrightarrow \lambda \mapsto (\Psi_\vartheta)^p(\lambda)=\lambda^{p\vartheta} \in \Lle (\Rset+)$$
    and the latter is equivalent to  $1\leq p\vartheta= p/\alpha\leq 2.$
\end{enumerate}
\end{remark}
\section{Some results on Voiculescu transforms related to Hyperbolic functions} \label{app}
In the following, for Voiculescu (Nevanlinna) transforms found in Jurek \cite{j2}, we compute the triplet of characteristics $(\aar_{\Xs}, \, \bbr_{\Xs},\,\Fou[\rho_{\Xs},.])$  in  their corresponding representations \eqref{(1)} via the inversion formula  \eqref{(31)}. For more  facts and  formulas, we refer to the Appendix at the end of this article. We recall that $\tilde{X}$ indicates the free-probability analog, in the sense of the procedure {\bf(P4)} in Section \ref{free}, of the classical hyperbolic  characteristic function $X=C,\, S, T$ which are self-decomposable, that
\begin{equation}\label{hyp}
\Four_{\Cs}(s)=\frac{1}{\cosh(s)}, \quad \Four_{\Ss}(s)=\frac{s}{\sinh(s)}, \quad \Four_{\Ts}(s)= \frac{\tanh(s)}{s}, \quad s\in \Rset.
\end{equation}
and that the Voiculescu transforms of propositions \ref{pro1}, \ref{pro2} and Corollary \ref{pro3} correspond to the free-probability analogs $\tilde{X}$.
\begin{proposition}\label{pro1} Recall the $\beta$ function is defined in \cite[\textbf{8.371(2)}]{grad} by
\begin{equation} \label{bet}
\beta(z)=\int_0^\infty e^{-z x}\,\frac{dx}{1+e^{-x}},\quad z\in\Cset^+.
\end{equation}
For the free-infinitely divisible Voiculescu transform
$$V_{\Cts}(iw)= -F_{\Cts}(iw)=i\left[1-w\beta\left(\frac{w}{2}\right)\right], \quad w>0,$$
the characteristics are given by   $\aar_{\Cts}=\bbr_{\Cts}=0$,  and the measure $\rho_{\Cts}$ is a such that $\rho_{\Cts}(\Rset)= \frac{\pi}{2} -1$  and
\begin{equation}\label{(5)}
\Fou[\rho_{\Cts}\,;\,s]=2\,\sinh(s)\,\arctan\,(e^{-s})+\frac{\pi}{2} e^{-s}-1=  \int_0^\infty\cos(sx) \,\frac{x}{(1+x^2)\sinh\left(\frac{\pi}{2} x\right)}\, dx,\quad  s \in \Rset,
\end{equation}
where $\beta$ is  given by \eqref{bet} in the Appendix.
\end{proposition}
\begin{proposition} \label{pro2} Recall the digamma function   is given by
\begin{equation} \psi(z)=\log z +\int_0^\infty\big(\frac{1}{s}-\frac{1}{1-e^{-s}}\big) e^{-  z s}ds, \quad z\in \Cset^+,\qquad \cite[\textbf{8.361(8)}]{grad}
\label{psii}\end{equation}
and the Euler-Mascheroni constant corresponds to $\gamma:=- \psi(1)$. For the free-infinitely divisible Voiculescu transform
$$V_{\Sts}(iw)=- F_{\Sts}(iw)= i\left[w\,\psi\left(\frac{w}{2}\right)-w\,\log\left(\frac{w}{2}\right) +1\right] ,\quad  w> 0,$$
the characteristics are given by  $\aar_{\Sts}=\bbr_{\Sts}=0$ and the measure $\rho_{\Sts}$ is  a such that $\rho_{\Sts}(\Rset)=\gamma + \log2-1$ and
\begin{eqnarray}
\Fou[\rho_{\Sts}\,;\,s] &=&  \frac{e^{-s}\Ei(s) +e^{s}  \Ei(-s)}{2}+    \cosh(s)\, \log \left(\frac{1+e^{-|s|}}{1-e^{-|s|}}\right)-1  \nonumber\\
&=&  2 \int_0^\infty \cos(sx) \frac{x}{(1+x^2)\left(e^{\pi x}-1\right)} \, dx, \quad s \in \Rset.
\label{62} \end{eqnarray}

where  $\Ei$ is the {\it exponential integral function} given by \eqref{ei} in the Appendix.
\end{proposition}

By the elementary relation $\Four_{\Cs} =\Four_{\Ss} \cdot \Four_{\Ts}$, we have $C\simdis S+T$ , where $S$ and $T$ independent versions. Using Remark \ref{simdis},  we can state:
\begin{corollary}\label{pro3}  For the free-infinitely divisible Voiculescu transform,
$$V_{\Tts}(iw)=-F_{\Tts}(iw)= F_{\Sts}(iw)-F_{\Cts}(iw)= iw\,\left[\,\log \left(\frac{w}{2}\right)-\beta\left(\frac{w}{2}\right)-\psi\left(\frac{w}{2}\right)\,\right],\quad  w>0,$$
the  characteristics are given by   $\aar_{\Tts}=\bbr_{\Tts}=0$,  and the measure $\rho_{\Tts}$ is  a such that $\rho_{\Tts}(\Rset)= \frac{\pi}{2} - \gamma -\log2$  and
\begin{eqnarray}
\Fou[\rho_{\Tts}\,;\,s]&=& \frac{\pi}{2} e^{-|s|}+ 2 \sinh(s)\arctan(e^{-|s|})+\cosh (s) \, \log\frac{1-e^{-|s|}}{1+e^{-|s|}} - \frac{[e^{-s} \Ei(s)+e^{s} \Ei(-s)]}{2}    \nonumber \\
&=&   \int_0^\infty \cos(sx) \frac{2x}{(1+x^2)\left(e^{\frac{\pi}{2} x}+1\right)} \, dx, \quad s\in \Rset\label{(7)}.
\end{eqnarray}
\end{corollary}
\medskip

In general, if $X$ is a self-decomposable random variable,  then the corresponding  L\'evy measure $\nu_{\Xs}$ has the form
$$\nu_{\Xs}(dx)= h_{\Xs}(x)\,  dx , \quad x\neq 0$$
$h$ is a measurable function such that  $x\mapsto x h_{\Xs}(x)$ is increasing on $(-\infty,0)$ and decreasing on $(0,\infty)$.
Therefore, there exists a BDLP $Y=Y_{\Xs}$, from the corresponding random integral representations \eqref{bdl}, such that the L\'evy  measures is $\nu_Y(dx)=-(x\,h_{\Xs})'(x)\,dx$ whenever $h_{\Xs}$ is differentiable, cf. Steutel and van Harn \cite[Proposition 6.11 and Theorem 6.12, Chapter V]{steutel} and also Jurek \cite[Corollary 1.1, p.97]{j1}, \cite[Section 2.1]{j2} or Jurek and Yor \cite[p. 183, formulae (d) and (e)]{jy}.
Consequently, on the level of  Nevanlinna measures, we  have, by \eqref{level}
\begin{equation}\label{back}
\rho_{\tilde{Y}_{\Xs}}(dx):=\frac{x^2}{1+x^2}\nu_Y(dx) = - \frac{x^2}{1+x^2}\,\big(h_{\Xs}(x)+x h_{\Xs}^\prime(x)\big) dx = - \rho_{\Xts}(dx)- \frac{x^3\,  h_{\Xs}^\prime(x)}{1+x^2} dx.
\end{equation}
where $\tilde{X}$ and $\tilde{Y}_{\Xs}$ are the free analogues of $X$ and $Y_{\Xs}$, respectively. As in the previous propositions we have a similar results for the BDLP's as well, although we computed  it  for $\tilde{Y}_{\Cs}$, only:
\begin{proposition}\label{pro4} Let $K$ stands for the Catalan constant, and $\zeta$, $Li_n$  be the {\it Riemann's zeta  function} and the {\it polylogarithm functions}   given in \eqref{rz} and \eqref{pl} in the Appendix. For the free-infinitely divisible Voiculescu transform
$$ V_{\tilde{Y}_{\Cs}} (iw)= -F_{\tilde{Y}_{\Cs}} (iw) = i\,\left[1+\frac{w^2}{2}\,\zeta\left(2, \frac{w}{2}\right)-\frac{w^2}{4}\,\zeta\left(2, \frac{w}{4}\right)\right], \quad w>0, $$
the characteristics are given by  $\aar_{\tilde{Y}_{\Cs}}=\bbr_{\tilde{Y}_{\Cs}}=0$ and  the measure $\rho_{\tilde{Y}_{\Cs}}$ is such that $\rho_{\tilde{Y}_{\Cs}}(\Rset)= 2\,K-1$  and
\begin{eqnarray}
\Fou[\rho_{\tilde{Y}_{\Cs}}\,;\,s]&=&2\,\cosh(s)\left(K-\int_0^{s}\frac{x}{\cosh(x)}dx\right)-s\tanh(s)-1 \nonumber\\
&=&\frac{\pi}{2} \int_0^\infty \cos(sx)\,\frac{x^2}{1+x^2}\,\frac{\cosh(\frac{\pi x}{2})}{\sinh^2( \frac{\pi x}{2})}dx, \quad s \in \Rset .
\label{fl2} \end{eqnarray}
\end{proposition}
\medskip

As a by-product of our Propositions \ref{pro1} and \ref{pro4}, we have
\begin{corollary}\label{cor1} With the notations of Proposition \ref{pro4}, we have
$$\Fou[\rho_{\tilde{Y}_{\Cs}}\,;\,s]+\Fou[\rho_{\Cts}\,;\,s]=2\int_0^\infty\cos(sx)\frac{x^3}{1+x^2}\big( -h_{\Cs}(x)\big)^\prime\, dx, \quad  s \in \Rset,$$
\noindent where the function $h_{\Cs}(x):= 1/\left( 2x \sinh (\frac{\pi}{2} x)\right)$ is the density of the L\'evy measure of the hyperbolic cosine  function  $\Four_{\Cs}$.
\end{corollary}
\begin{remark}  Formula \eqref{(5)} is confirmed by \cite[\textbf{4.113}(8)]{grad} and \eqref{62}  was confirmed numerically for $s= 0.5,\, 1,\,2$. Formulae \eqref{(7)} and \eqref{fl2} seem  to be new and  might be of some  interest.
\end{remark}
\section{Proofs}\label{proo}
\subsection{Useful results on cumulant functions} \label{useful}
Consider   the class of cumulant functions denoted by
$$\CF =\{\lambda \mapsto \phi_{\Xs}(\lambda)= -\log\er[ e^{-\lambda X}],\; \, \lambda \geq 0, \quad  \mbox{where $X\,$ is a non-negative r.v.} \}.$$
\begin{remark}\label{discu}
The class $\CF$ contains the class of Bernstein functions $\BF$. By injectivity of the Laplace transform, it is seen that to every $\phi\in \CF$ corresponds a unique  non-negative r.v. $X$ such that $\phi=\phi_{\Xs}$, i.e.
\begin{equation}
\phi(\lambda) = -\log\er[ e^{-\lambda X}]= -\log \big(\pr(X=0) +  \pr(X>0)\;\er [e^{-\lambda X}\mid X> 0]\big), \quad \lambda \geq 0.
\end{equation}
Observe that $\phi$ is linear if, and only if, $X$ is deterministic. Also observe that

\begin{itemize}
\item $\phi$ is infinitely  differentiable $(0, \infty)$ and is a strictly increasing bijection $[0, \infty)\to [0,  l_\phi)$. We denote from now on, by
    $$\Psi:[0,l_\phi)\to \Rset_+, \; \,\mbox{the inverse function of} \;\, \phi; $$
\item we necessarily have
\begin{equation}\label{fifi}
l_\phi := \lim_{\lambda \to \infty}{\phi(\lambda)}=- \log\pr(X=0);
\end{equation}
\item if $\,\Ebf\,$ denotes a r.v. with standard exponential distribution, then
$$\pr\left(\frac{\Ebf}{X}>x\right)=\pr(X=0)+\pr(X>0)\;\er[e^{-xX}\mid X>0],\quad\mbox{for all $x\geq 0$}$$
and
$$\pr\left(\Psi \big( \min (\Ebf, l_\phi)  \big) >x \right) =\pr \left(  \min (\Ebf, l_\phi)  >\phi(x) \right)= \pr (\Ebf > l_\phi )+ \pr \big(\Ebf < l_\phi,  \Ebf  >\phi(x) \big)= e^{-\phi(x)}.$$
\end{itemize}
\end{remark}
\medskip
Exploiting last remarks, we propose the following lemma that will be useful in the sequel:
\begin{lemma}[Some results on cumulant functions] \label{keylemma} Let  $\phi:\Rset_+ \to \Rset_+$ a differentiable function on $(0,\infty)$, $\phi(0)=0$. Denote $l_\phi := \lim_{\lambda \to \infty}{\phi(\lambda)}$ and by $\,\Ebf\,$  a r.v. with standard exponential distribution.
\begin{enumerate}[1)]
\item  Assume $\phi$ is a non-linear cumulant function, hence is  associated to a  non-negative and non-deterministic r.v. $X$. Then,
      \begin{enumerate}[(i)]
      \item   $\phi'$ is a (strictly) decreasing bijection on $(0,\infty)$, hence $\phi$ is  strictly concave and increases to $l_\phi$;
      \item The inverse function $\Psi$ of $\phi$, is such that $\Psi':(0,l_\phi) \to (0,\infty)$ is a (strictly) increasing bijection, hence $\Psi$ is strictly convex;
      \item As $\lambda \uparrow \infty$, the function  $\phi(\lambda)/\lambda$   decreases to
        $$L_X:=\inf\{x\geq 0, \; \mbox{s.t.}\; \pr(X\leq x)>0\}.$$
        Hence, on $(0,\infty)$, the function $\Psi(x)/x$ increases to $1/L_X$ as $x \uparrow l_\phi$.
        Further,
        \begin{equation}  \label{debloc11}
         L_X=\lim_{x\to \infty}\phi'(x) \in [0,\infty);
         \end{equation}
      \item If $p >0$, then $(\Psi^p)'$ is positive on $(0,\infty)$ and we have the equivalence
          \begin{equation}\label{ppp}
          (\Psi^p)'(0+):= \lim_{\lambda \to 0+}   (\Psi^p)'(\lambda)<\infty \Longleftrightarrow \lim_{\lambda \to 0+} \frac{\lambda^{p-1}}{\phi'(\lambda)};
           \end{equation}
      \item The function  $\Psi^p$, is    convex (respectively concave) if, and only if,
         $$\lambda \mapsto \lambda^{1-p} \,\phi' \; \mbox{is non-increasing (respectively non-decreasing )}.$$
         In  particular, if $p\geq 1$, then  $\Psi^p$ is  strictly convex  and $(\Psi^p)'(0+)\in [0,\infty)$.
      \end{enumerate}
\item Assume $\phi$ is a strictly  increasing  bijection with inverse  $\Psi$.  Then, the following statements are equivalent.
    \begin{enumerate}[i)]
    \item $\phi$ is a cumulant function;
    \item For some non-negative r.v. $X$, we have the identity in law
       \begin{equation} \frac{\Ebf}{X} \simdis \Psi \left(\min \left(\Ebf, l_\phi\right)\right),\qquad \mbox{(where $\Ebf$ in the l.h.s. is assumed to be independent of $X$};
      \label{id1}\end{equation}
    \item For some non-negative r.v. $X$, we have the representations
       \begin{equation}\label{integ0}
       \er[X^{-p}] =  \frac{M(p)}{\Gamma(p+1)}, \qquad M(p):= \int_0^{l_\phi} e^{- x} \,\Psi^p(x) \,dx + e^{-l_\phi}\, \Psi^p(l_\phi)\, \II_{(l_\phi <\infty)},
       \end{equation}
        for every $p  >-1,\,$ such that $\, e^{- x}\,\Psi^p(x)\,$ is integrable on $(0,l_\phi)$.
    \end{enumerate}
\end{enumerate}
\end{lemma}
\begin{remark} \label{remel}  Observe that formula \eqref{integ0} reminds the one involving the Gamma function in \eqref{gag}. Actually, under the condition of  integrability  of $\Psi^p$,   the injectivity of the Mellin transform insures that  the function $M(p), \; p>-1$ is the Mellin transform  of a (unique) positive random variable, namely
   \begin{equation}\label{mellin} M(p)=  \er \left[ \Psi^p\big(\min (\Ebf, l_\phi)\big)\right] ,
   \end{equation}
   thus, \eqref{integ0} reads
   $$ \er \left[\,\left(\frac{\Ebf}{X}\right)^p\,\right]=  \er[\Ebf^p]\; \er[X^{-p}]= \Gamma(p+1)\; \er[X^{-p}]= \er \left[ \Psi^p\big(\min (\Ebf, l_\phi)\big)\right].$$
\end{remark}
\begin{proofs}{\bf of Lemma \ref{keylemma}.}  1) $(i)$:   The first and the second derivative of $\phi$ are given by
\begin{equation}\label{p2}
\phi'(\lambda)=  \frac{\er[X e^{-\lambda X}]}{\er[ e^{-\lambda X}]} \quad \mbox{and }\quad \phi''(\lambda) = \frac{\er[X e^{-\lambda X}]^2- \er[X^2 e^{-\lambda X}]\er[e^{-\lambda
X}]}{\er[ e^{-\lambda X}]^2}, \qquad  \lambda >0.
\end{equation}
By Cauchy-Schwarz inequality, every pair of nonnegative random variable $Y,Z$ satisfy the inequality  $\er[yZ]^2\leq \er[Y^2] \er[Z^2]$ and the equality holds if and only if there exists $c\geq0$ such that $Y= c Z$ a.s. Choosing  $Y=X e^{-\frac{\lambda}{2} X}$ and $Z= e^{-\frac{\lambda}{2} X}$ in \eqref{p2}, we recover the negativity of   $\phi''$ and the claim on $\phi'$ follows.\\
1) $(ii)$: The assertion is due to the observation $\Psi'=1/\phi'(\Psi)$ and $\Psi''=-\phi''(\Psi)/(\phi'(\Psi))^3$.\\
1) $(iii)$ The assertion is an adaptation of the result of Pakes \cite[Theorem 2.1]{pakesc} by taking the Mellin transform   $M(\lambda)^{1/\lambda}=\er[(e^{-X})^\lambda]^{1/\lambda}=e^{-\phi(\lambda)/\lambda}$ there. The limit in  \eqref{debloc11} is obtained by Karamata's Theorem  \cite[Theorem 1.5.11]{BIN}, the monotone density theorem \cite[Theorem 1.7.2]{BIN}  and the concavity of a cumulant function.\\
1) $(iv)$: We have
\begin{equation}\label{pp}
   (\Psi^p)' = p \;\frac{\Psi^{p-1} }{\phi'\left(\Psi \right)}, \qquad  (\Psi^p)''= \frac{p \;\Psi^{p-2}}{\phi'(\Psi)^3} \left[(p-1)\,\phi'(\Psi)-\Psi  \phi''(\Psi)\right]
  \end{equation}
Positivity of $ (\Psi^p)'$ is immediate and since $\Psi(0)=0$,  then by the change of variable $x=\phi(\lambda)$, we have
\begin{equation}\label{pp0}
\lim_{x\to 0+}(\Psi^p)'(x)= p\; \lim_{x\to 0+} \frac{\Psi^{p-1}(x)}{\phi'(\Psi(x))}=   p\; \lim_{\lambda\to 0+} \frac{\lambda^{p-1} }{\phi'(\lambda)}.
\end{equation}
For a later use, one can observe  that
\begin{equation}\label{lpp}
p\geq 1  \Longrightarrow \lim_{x\to 0+}x\;(\Psi^p)'(x) = p \;\lim_{\lambda \to 0+} \phi(\lambda) \frac{ \lambda^{p-1}}{\phi'(\lambda)} \in [0,\infty).
\end{equation}
1) $(v)$:   By \eqref{pp}, we have
$$Sign\left((\Psi^p)''\right)= Sign\left[(p-1)\,\phi'-\lambda  \phi'' \right]= Sign\left[\left(  \frac{\lambda^{p-1}}{\phi'}\right)'\right] $$
and if $p\geq 1$, then  $\lambda^{p-1}/\phi'$ is increasing, viz.  $\Psi^p$ is convex and last assertion is then obtained with the help of \eqref{pp0}. \\
\noindent 2) $(i)\Longleftrightarrow (ii)$: Is obtained from the discussion in Remark \ref{discu}.\\
2) $(ii) \Longrightarrow (iii)$: The  case $l_\phi=\infty$ is easy to prove due to the  facts, that $\Psi(0)= \phi(0)=0$,  $\Psi$ is an increasing bijection from $\Rset_+\to \Rset_+$ (with inverse $\phi$), and that
$$ - \frac{d}{dy}\er[e^{-y X}] = - \frac{d}{dy}e^{- \phi(y)}=  \phi'(y)\, e^{-  \phi(y)}= \er[X\, e^{-y X}].$$
Indeed, making the change of variable $x=\phi(y)$, using Tonelli-Fubini's theorem  and the representation of the gamma function \eqref{gag}, these facts entail
\begin{eqnarray*}
 \int_0^\infty e^{-x}  \,\frac{\Psi^p(x)}{\Gamma(p+1)} \,dx&=&  \int_0^\infty  \phi'(y) \,e^{-  \phi(y)} \,  \frac{y^p}{\Gamma(p+1)} \,dy  =  \int_0^\infty \er[X\, e^{-y X}] \,  \frac{y^p}{\Gamma(p+1)} \,dy \\
&=& \er\left[\frac{1}{X^p}  \underbrace{ \int_0^\infty \, e^{-y X} \,  \frac{X^{p+1}\,y^p}{\Gamma(p+1)} \,dy}_{=1} \right] = \er\left[X^{-p} \right].
\end{eqnarray*}
The  case $l_\phi<\infty$ is  proved similarly.\\
2) $(iii)\Longrightarrow(ii)$: This has been noticed in Remark \ref{remel}.
\end{proofs}
\subsection{The proofs}
\begin{proofs}{\bf of Theorem \ref{vraiz}.} Observe that if $F$ is a Nevanlinna-Pick function represented by \eqref{n1}, then  $\aar, \, \bbr$ and $\rho(\Rset)$ are given by \eqref{recup}.\\
\noindent  1) $\Longrightarrow$ 2): It is enough to reproduce the steps of the proof \cite[Corollary 6.13]{SSV} taking into account that $\rho$ is supported by $A$.\\
\noindent 2) $\Longrightarrow$ 3):  Use Jankowski and Jurek inversion's procedure given \cite[Theorem 1]{jan} and explained before \eqref{(1)}  and observe that  $H(w)=G(iw), \; w>0$. Also observe that
\begin{equation}\label{(4)}
\mathfrak{L}\big[\sinh (x); w]= \frac{1}{w^2-1}  \quad  \mbox{and}  \quad   \mathfrak{L}[\cosh (x); w]=\frac{w}{w^2-1},  \quad  w>1,
\end{equation}
which yields  \eqref{(31)}.\\
\noindent 3) $\Longrightarrow$ 1): Observe that knowing the function $w\mapsto F(iw)=i(1-w^2)\, H(w)$ on the interval $(1, \infty)$ or even on any interval $(a,b)\subset (0,\infty)$ is sufficient to fully recover $F$, cf. Jedidi et al. \cite[Lemma 4.1]{JBH}.

\smallskip
\noindent The last claim is a straightforward adaptation  of \cite[Corollary 6.13]{SSV}.
\end{proofs}
\medskip

\begin{proofs}{\bf of Proposition \ref{zfz}.} In order to get \eqref{fz}, just use the representation \eqref{LE} of $\Psi$,   the properties  \eqref{gag} of the Gamma function and the identity
$$-z^2  \left(\frac{1}{z+x}-\frac{1}{z}+\frac{1}{z^2}\frac{x}{1+x^2}\right)=\frac{z x -1}{z+x}\, \frac{x^2}{1+x^2},$$
in order to write: for $z\in \Cset^+$
\begin{eqnarray*}
F(z)&=&z^2\int_0^\infty e^{-z u}\,\left(\bbr  -  \beta u  -\frac{\gamma^2}{2} u^2 - \int_{\oi}\left(e^{-u x}-1+ u \frac{x}{1+x^2}\right)\nu(dx)\right) \,du\\
&=&-\beta +\bbr z-\frac{\gamma^2}{z}-z^2 \int_{\oi} \left(\frac{1}{z+x}-\frac{1}{z}+ \frac{1}{z^2} \frac{x}{1+x^2}\right)\nu(dx)\\
&=&-\beta +\bbr z+ \int_{[0,\infty)} \frac{z x -1}{z+x}\, \left( \gamma^2 \delta_0(dx)+ \frac{x}{1+x^2} \nu(dx)\right)\\
&=&  \aar +\bbr z+\int_{[0,\infty)} \frac{z x -1}{z+x}\rho(dx).
\end{eqnarray*}
The controls \eqref{cont}  and Tonelli-Fubini's theorem allow the reversal of order in the integrals.
\end{proofs}
\medskip

\medskip

\begin{proofs}{\bf of Proposition \ref{prop-moment}.} 1) $\Longrightarrow$ 2): If $\phi \in \BF$, then it is associated to a subordinator $\xi={(\xi_t)}_{t\geq 0}$, such that each $\xi_t,\;t>0$, has the cumulant (Bernstein) function $\phi_t=t\,\phi$ whose inverse is $\Psi_t(x)=\Psi(x/t)$. Since $\phi(\infty)=\infty$, then assertion  \eqref{id2} immediately follows from \eqref{id1}.\\
2) $\Longrightarrow$ 3): is trivially satisfied with $\eta=\xi$.\\
3) $\Longrightarrow$ 1): By 2) in Lemma \ref{keylemma}, necessarily the inverse $t \phi$ of $\Psi_t(x)=\Psi(x/t)$ is the cumulant of $\eta_t$ for all $t>0$, i.e., $e^{-t\phi}$ is a completely monotone, or equivalently
$(1-e^{-t\phi})/t$ is a Bernstein function for all $t>0$. Passing to the limit as $t\to 0$ and using the fact that the class $\BF$ is closed under pointwise limits, cf. \cite[Corollary 3.8]{SSV}, we deduce that $\phi \in \BF$.\\
3) $\Longrightarrow$ 4): By  the  implication 2) $\Longrightarrow$ 3), point 3) is equivalent to the existence of the subordinator $\xi$ satisfying \eqref{id2}, then \eqref{formulei} follows from of \eqref{integ0} and from the change of variable $x\rightarrow  tx$ there.\\
4) $\Longrightarrow$ 1): By Lemma \ref{keylemma}, necessarily the inverse $t \phi$ of $\Psi_t(x)=\Psi(x/t)$ is the cumulant of $\xi_t$ and  $\xi$ is a subordinator. Then, $\phi$ is a Bernstein function.
\end{proofs}
\medskip

\begin{proofs}{\bf of Theorem \ref{temporalcm}.} 1) Property 1) $(iii)$ in Lemma \ref{keylemma}  gives the implication
$$x\mapsto \frac{\Psi(x)}{x} \;\mbox{is non-decreasing}, \quad p\in(-1,0)\Longrightarrow e^{-tx}\Psi^p(x) \;\mbox{is integrable at $\infty$,}\quad \forall t>0.$$
Then, the definiteness of  $\Lap[\Psi^p\,;\,.]$  on $\oi$   is equivalent to   1)$(i)$. Since $\Psi^p$ is a non-increasing function, the latter is equivalent to assertion 1)$(ii)$, because
representation \eqref{formulei} of $\er[\xi_t^{-p}]$ meets the form \cite[(3.3) p. 23]{SSV}  of a Bernstein function.  Representation \eqref{formula1}  is obtained by integration by parts fom \eqref{formulei}.\\
2) Similarly,
$$x\mapsto \frac{\Psi(x)}{x} \;\mbox{is non-decreasing}, \quad p \geq 0\Longrightarrow \Psi^p(x) \;\mbox{is integrable at}\;0.$$
Then, the definiteness of  $\Lap[\Psi^p\,;\,.]$  $\oi$ is then  equivalent to 2)$(i)$.
Since $\Psi^p$ is an non-decreasing function, null at 0, and by an integration by parts, the latter is equivalent to representation \eqref{cmety}, hence $\er[\xi_t^{-p}]$ meets the form \label{cm1}  of a completely monotone function.
The last claim  is due  to  1) $(v)$ in Lemma \ref{keylemma} that insures the concavity of $\Psi^p$, and to the form \eqref{formulei} that gives
\begin{equation}\label{forp}
t \;\er[\xi_t^{-p}] = \frac{t^2}{\Gamma(p+1)}\int_0^\infty e^{-t x} \; \Psi^p(x) \,dx,
\end{equation}
which meets the form  \cite[(3.4) p. 23]{SSV} of a Bernstein function. Representation \eqref{formula2} is obtained by integration by parts in \eqref{cmety}. \\
3) The convexity of $\Psi^p$ and $(\Psi^p)'(0+)\in [0,\infty)$ are guaranteed by   1)$(v)$  in Lemma \ref{keylemma}. The integrability of $\Psi^p$ near zero is guaranteed by $p\geq 1$. If further $\Lap[\Psi^p\,;\,t]<\infty$, for some $t>0$, then   two integration by parts in \eqref{forp}, give representation
\begin{equation}\label{formulei2}
t\; \er[\xi_t^{-p}] \;  = \frac{1}{\Gamma(p+1)} \left((\Psi^p)'(0)+ \int_{0}^\infty e^{-t x} \, (\Psi^p)''(x) \,dx \right)
\end{equation}
Differentiation
$$-\frac{d}{dt}\left(t\; \er[\xi_t^{-p}] \right)= \frac{1}{\Gamma(p+1)}\int_0^\infty e^{-t x} \;x \; (\Psi^p)''(x) \,dx.$$
The function $t \;\er[\xi_t^{-p}]$ being positive, deduce the equivalences $3)(i) \Longleftrightarrow  3)(ii)$ from $2)(i) \Longleftrightarrow 2)(ii)$.

 The last claims in 3) are obtained as follows:\\
3) $(a)$ Representation   \eqref{cmety}, when multiplied by $t^2$  meets with the form  \cite[(3.4) p. 23]{SSV} of a Bernstein function, if, and only if $(\Psi^p)'$ is concave.\\
3) $(b)$ Representation \eqref{repp4} is a simple reformulation of  \eqref{rerp1}, when differentiating  twice  the expression of $\Psi^p$ in \eqref{LE}.
\end{proofs}
\medskip

\begin{proofs}{\bf of Corollary \ref{tempo}.} 1)$(i)$ and  1)$(ii)$ are a direct consequence of Proposition \ref{prop-moment}  because $\xi$ is a subordinator.\\
1)$(iii)$:  If $\Psi^p \in \Lle(0,\infty)$, then, by Proposition \ref{prop-moment}, the function $f(t) :=   t \, \er[\xi_t^{-p}], \; t >0,$ has the representation
$$f(t)=\frac{t^2}{\Gamma(p+1)} \int_0^\infty e^{-t x} \,\Psi^p(x)\, dx.$$
By Proposition \ref{zfz}, the  extension on $\Hset^+$ of  $f$
is such that $-f$ is a non-positive Nevanlinna $\Nev(\Rset_+)$-function, which is continuous on $\oi$. By \cite[Corollary 6.14]{SSV} this equivalent to say that $f$ is a Stieltjes transform and by \cite[ Theorem 6.2]{SSV}, this also equivalent to  $t \mapsto t \, f(t) \in \CBF$.\\
\medskip
\noindent 2)  Since $\xi$ is a subordinator, then, by Proposition \ref{prop-moment},  the function $t \mapsto f(t)=:t\; \er[\xi_t^{-p}]$ is represented by
$$ -f(t) = - t^2 \int_0^\infty e^{-t x} \,\Psi^p(x)\, dx.$$
Since $f$ is a Stieltjes function, then, using \cite[Corollary 6.14]{SSV} again, we deduce the extension on $\Hset^+$ of $-f$ is a Nevanlinna $\Nev^\star(\Rset_+)$-function because it satisfies
$\lim_{|z|\to \infty, z\in \Cset^+  }-f(z)/z =0$. Finally, by Proposition \ref{zfz}, we necessarily have $\Psi^p \in \Lle(0,\infty)$.  To conclude, observe that the inverse
function of $\Psi^p$ is the Bernstein function $(\Psi^p)^{-1}(\lambda)=\phi(\lambda^{1/p})$ corresponding to the subordinator $\xi^{(\scriptstyle 1/p)}\circ  \xi$.
\end{proofs}
\medskip

In all the following proofs, the symbol $\ast$ denotes the usual additive convolution of functions on the positive half-line. The calculi where also checked by WolframAlpha or Mathematica.

\medskip
\begin{proofs}{\bf of Proposition \ref{pro1}.}
Firstly,  note that $F_{\Cts}(i)=-i \big(1-\beta(1/2)\big)= i(\pi/2-1)$ and therefore  $\rho_{\Cts}(\Rset)=\frac{\pi}{2}-1$. Secondly, by \eqref{bet}, we have $\beta(z)=  \mathfrak{L}[ (1+e^{-x})^{-1}; z], \; z \in \Cset^+.$ Consequently,
\begin{eqnarray}
\frac{iF_{\Cts}(iw)}{w^2-1}&=&\frac{1- w\beta(\frac{w}{2})}{w^2-1}=\mathfrak{L}[\sinh (s); w] - \mathfrak{L}[\cosh (s)\,;w\,]\;\mathfrak{L}\left[\frac{2}{1+e^{-2s}}\,;w\,\right]\nonumber \\
&=&\mathfrak{L}[\sinh (s); w] -\mathfrak{L}\left[\, \left(\cosh(u)\ast \frac{2}{1+e^{-2u}}\right)(s); w\right],\quad \mbox{for}\; \;w>0, \, w\neq 1.\label{(9)}
\end{eqnarray}
Thirdly,  by a differentiation, one checks that, for $s>0$,
\begin{equation*}
\left( \cosh(u)\ast  \frac{2}{1+e^{-2u}}  \right)(s) := \int_0^s\cosh(s-u)\,\frac{2}{1+e^{-2(s-u)}}du =2\sinh(s)\, \left[\frac{\pi}{4} - \arctan(e^{-s})\right]-e^{-s}+1,
\end{equation*}
and inserting the latter into \eqref{(9)}, we get
\begin{equation}
\frac{iF_{\Cts}(iw)}{w^2-1}=\mathfrak{L}\left[\sinh (s)-2\sinh (s) \left(\frac{\pi}{4}-\arctan(e^{-s})\right)+e^{-s}-1\,;w\,\right].
\end{equation}
Finally, since $\rho_{\Cts}(\Rset)=\frac{\pi}{2}-1$ and using  \eqref{(4)}, we arrive to
\begin{eqnarray*}
\Fou[\rho_{\Cts}\,;\,s] &=& \left(\frac{\pi}{2}-1\right)\cosh (s)+ \sinh (s) -2\sinh (s) \,\left(\frac{\pi}{4} -\arctan(e^{-s})\right)+e^{-s}-1 \\
&=&  2\,\sinh (s) \,\,\arctan(e^{-s})+\frac{\pi}{2}  e^{-s}-1,
\end{eqnarray*}
which gives the first equality in  \eqref{(5)}.  On the other hand, from Jurek \cite[Corollary 2]{j2},  we know that $V_{\Cts}(iw)=-F_{\Cts}(iw)$ is a free-probability analog of the classical hyperbolic characteristic function $1/\cosh(t)$ whose Nevanlinna measure $\rho_{\Cts}$ has a density
$$ \frac{|x|}{2(1+x^2)\sinh(\frac{\pi}{2}|x|)}, \quad  x\in \Rset.$$
Thus,
$$ \Fou[\rho_{\Cts}\,;\,s]=\int_{\Rset} e^{isx} \frac{|x|}{2(1+x^2)\sinh(\frac{\pi}{2}|x|)} \,dx= \int_0^\infty \cos (sx)\frac{x}{(1+x^2)\sinh(\frac{\pi}{2}x)} \,dx,$$
which shows the second equality in  \eqref{(5)}.
\end{proofs}
\medskip

\begin{proofs}{\bf of Proposition \ref{pro2}.}
By \cite[\textbf{8.366(2)}]{grad}, we have $\psi(1/2)= - \gamma -2 \log2$,  then $F_{\Sts}(i)= i(\gamma +\log 2-1)$.  Hence, in \eqref{(1)}, we have $\aar_{\Sts}=\bbr_{\Sts}=0$ and the  measure $\rho_{\Sts}$  has a finite mass   $\rho_{\Sts}(\Rset)= \gamma + \log 2-1$.  Using the  integral formula  \eqref{psii} for the $\psi$-function, we get
\begin{eqnarray}
\frac{i F_{\Sts}(iw)}{w^2 -1}&=&\frac{1+w \,(\psi(w/2)-\log(w/2))}{w^2-1}=  \mathfrak{L}[\sinh (s); w] +\mathfrak{L}[\cosh (s); w] \;\mathfrak{L}\left[ \frac{1}{s}-\frac{2}{1-e^{-2s}}; w\right] \nonumber\\
&=& \mathfrak{L}\left[\sinh (s) + g(s); w\right],\quad \mbox{where} \;g(s):=\left(\cosh(u) \ast \big(\frac{1}{u}-\frac{2}{1-e^{-2u}}\big)\right)(s).\label{(11)}
\end{eqnarray}
By a direct differentiation in identities (F3) in the Appendix of the $\Ei$ function, one checks that,  for $s>0$, we have
\begin{eqnarray*}
g (s) &=&\int_0^s \cosh(s-u) \left[\frac{1}{u} -\frac{2}{1-e^{-2u}}\right] \,du = \left[ \frac{e^{-s}}{2}\Ei(u)+\frac{e^s}{2}\Ei(-u) - e^{u-s}+\cosh(s)\, \log\frac{1+e^{-u}}{1-e^{-u}} \right]\Big|^{u=s}_{u=0^+}  \\
&=& \frac{e^{-s}}{2}\Ei(s)+\frac{e^s}{2}\Ei(-s) - 1+\cosh(s)\, \log\frac{1+e^{-s}}{1-e^{-s}} - \lim_{u \to 0^+} \Big[ \frac{e^{-s}}{2}(\Ei(u) -\log(1-e^{-u})) \\
&& \quad  +  \frac{e^s}{2}(\Ei(-u)-\log(1-e^{-u}))- e^{u-s}+ \cosh (u) \log(1+e^{-u})\big]\\
&=& \frac{e^{-s}}{2}\Ei(s)+\frac{e^s}{2}\Ei(-s) - 1 +e^{-s} + \cosh(s) \left( \log\frac{1+e^{-s}}{1-e^{-s}} -\gamma-  \log 2\right).
\end{eqnarray*}
Inserting the above equality into \eqref{(11)} and  using   \eqref{(4)} with $\rho_{\Sts}(\Rset)= \gamma + \log 2-1$, we get for $s\in \Rset$,
\begin{eqnarray*}
\Fou[\rho_{\Sts}\,;\,s]&=&(\gamma +\log2-1)\cosh(s)+\sinh(s) +\frac{e^{-s}}{2}\Ei(s)+\frac{e^s}{2}\Ei(-s)-1+e^{-s} +\cosh(s)\left(\log\frac{1+e^{-|s|}}{1-e^{-|s|}} -  \gamma - \log 2  \right)\\
&=&  \frac{e^{-s}}{2}\Ei(s)+\frac{e^s}{2}\Ei(-s) - 1 + \cosh(s) \, \log \frac{1+e^{-|s|}}{1-e^{-|s|}} ,
\end{eqnarray*}
which proves the  equality \eqref{62}.  From Jurek \cite[Corollary 3]{j2}, we know that $V_{\Sts}(iw)=-F_{\Sts}(iw)$ is the  free-analog of the classical hyperbolic sine characteristic function $\Four_{\Ss}(t)= t/\sinh(t)$ whose Nevanlinna measure equals to
$$\rho_{\Sts}(dx)= \frac{1}{2}\frac{|x|}{1+x^2} \frac{e^{-\frac{\pi}{2}|x|}}{\sinh(\frac{\pi}{2}|x|)}dx= \frac{|x|}{1+x^2}\, \frac{1}{e^{\pi |x|}-1}dx, \quad  x\in  \Rset, $$
and we get the  equality \eqref{62}.
\end{proofs}
\medskip

\begin{proofs}{\bf of Proposition \ref{pro4}.}  Since  $F_{\tilde{Y}_{\Cs}}(i)=i(2K-1)$,  then  in \eqref{(1)}, $\aar_{\tilde{Y}_{\Cs}}=\bbr_{\tilde{Y}_{\Cs}}=0$,  and for the measure $\rho_{\tilde{Y}_{\Cs}}$, we have  $\rho_{\tilde{Y}_{\Cs}}(\Rset)= 2K-1$. Using   \eqref{(4)} and  the integral representation \eqref{(8)} for the $\zeta(2, s)$  function,
we have
\begin{eqnarray*}
\mathfrak{L}[\Fou[\rho_{\tilde{Y}_{\Cs}}] (s)-(2K-1)\cosh(s)\,;w\,]&=&\frac{1}{w^2-1}\, \left[\frac{w^2}{2}\,\left(\zeta(2,\frac{w}{2})-
\frac{1}{2}\, \zeta(2,\frac{w}{4})\right)+1 \right]\\
&=& \frac{1}{w^2-1} +   2\frac{w^2}{w^2-1} \left[\int_0^\infty e^{-wu} \,\frac{u}{1-e^{-2u}}\, du -2 \int_0^\infty e^{-w u} \,\frac{u}{1-e^{-4u}}\, du \right]\\
&=&   \mathfrak{L}[\sinh (s) ; w] + 2\left(1+\frac{1}{w^2-1}\right)\,\int_0^\infty e^{-w u}\, u\;\frac{e^{-2u}-1}{1-e^{-4u}}\,du\\
&=& \mathfrak{L}[\sinh (s) ; w]+ (1+\frac{1}{w^2-1})\int_0^\infty e^{- wu}\, u\;\frac{1- e^{2u}}{\sinh(2u)} \,du \\
&=& \mathfrak{L}\left[\sinh (s)  + s\frac{1- e^{2s}}{\sinh(2s)}  ; w\right] +\mathfrak{L}[\sinh (s)\,;w\,]\;  \mathfrak{L}\left[s \frac{1-e^{2s}}{\sinh (2s)}\,;w\,\right]\\
&=& \mathfrak{L}\left[\sinh (s) + s \frac{1- e^{2s}}{\sinh(2s)} +h(s) ; w\right],
\end{eqnarray*}
where,
\begin{equation}\label{(13)}
h(s):= (\sinh u \ast u\frac{1-e^{2u}}{\sinh (2u)})(s)= \int_0^s\sinh(s-u)\,u\; \frac{1-e^{2u}}{\sinh (2u)}\,du.
\end{equation}
Using elementary computations, together with \eqref{latg}, \eqref{lig} in the Appendix, one checks that
\begin{eqnarray*}
h(s)&=& \left[e^{-s} (u-1)\,e^u  - i \cosh(s)\, \left(\Li_2(i e^u) - \Li_2(-i e^u)+ u\, \log \frac{1-ie^u}{1+ie^u}  \right)\right]_{u=0}^{u=s}\\
&=& \left[e^{-s} (u-1)\,e^u  -i \cosh(s)\, \left(-2i \int_0^u \frac{x}{\cosh(x)}dx \right)\right]_{u=0}^{u=s}\\
&=& e^{-s}  - 1+s  - 2 \cosh(s) \,\int_0^{s}\frac{x}{\cosh(x)}dx.
\end{eqnarray*}
\noindent From the above calculation and from the equality $(e^{2s}-1)/\sinh(2s)=2e^{2s}/(e^{2s}+1)$, we arrive at the representation of $\Fou[\rho_{\tilde{Y}_{\Cs}}], \; s\ge 0$:
\begin{eqnarray*}
\Fou[\rho_{\tilde{Y}_{\Cs}}] (s) &:=&(2K-1)\cosh(s) +\sinh(s) - 2s\,\frac{e^{2s}}{e^{2s}+1}+h(s)\\
&=&       2\,\cosh(s) \left(K-\int_0^{s}\frac{x}{\cosh(x)}dx\right)+ s\left(1- 2\frac{e^{2s}}{e^{2s}+1}\right) + \overbrace{\big(\sinh (s)-\cosh(s)+e^{-s}\big)}^{=0} -1\\
&=& 2\,\cosh(s) \left(K-\int_0^{s}\frac{x}{\cosh(x)}dx\right) -s  \tanh(s)-1\\
\end{eqnarray*}
which gives the equality \eqref{fl2}.  For the equality \eqref{fl2}, we just need to recall that  from Jurek and Yor \cite[Corollary 1 and a formula (7) p. 186]{jy}, that the Nevanlinna measure corresponding to the BDLP $Y_{\Cs}$, given by \eqref{bdl} is represented by
$$ \rho_{\tilde{Y}_{\Cs}}(dx)=\frac{ \pi}{4}\, \frac{x^2}{1+x^2}\,\frac{\cosh \left(\frac{\pi x}{2}\right)}{\sinh^2\left(\frac{\pi x}{2}\right)} dx, \quad x \in \Rset.$$
\end{proofs}
\medskip

\begin{proofs}{\bf of Corollary \ref{cor1}.} It is a simple application of Remark \eqref{back} applied to $X=C$.
\end{proofs}
\begin{remark}
\emph{Since the  hyperbolic sine and the hyperbolic tangent are self-decomposable as well, we may have a statement about $\tilde{S}$ and $\tilde{T}$, analogous to the one in Corollary \ref{cor1}  for the hyperbolic cosine function.}
\end{remark}
\section{Appendix}\label{appen}
\noindent For a convenience  of reading, we collect some facts on special functions. Boldface numbers below refer to formulae from \cite{grad}.
\begin{enumerate}[{\bf (F1)}]
\item Many of  those  functions we use are derived from \emph{Euler's $\Gamma$ gamma function} and the \emph{digamma function} $\psi$ used in Proposition \ref{pro2}:
    \begin{equation}\label{gag}
    \frac{\Gamma(p)}{z^{p}}:=\,\int_0^\infty e^{-zx} \,x^{p-1}\,dx,  \quad   \psi(p):= \frac{\Gamma'(p)}{\Gamma(p)}\quad z,\, p\in \Cset^+.
    \end{equation}
\item In Proposition \ref{pro1}, we used the  $\beta$-function
which was originally  derived from the digamma function via the formula
$$ \beta(x):=\frac{1}{2}\,\left[\psi\left(\frac{x+1}{2}\right)-\psi\left(\frac{x}{2}\right)\right] =\sum_{k=0}^{\infty}\frac{(-1)^k}{x+k},  \quad -x \notin \mathbb{N}  \qquad   \cite[\textbf{8.372(1)}]{grad}.$$

\item To formulate Proposition \ref{pro2}, we need  another special function, namely,  \emph{the exponential integral function} $\Ei$,  defined by
\begin{equation} \label{ei}
 \Ei(x):=-\int_{-x}^\infty\frac{e^{-t}}{t}dt, \quad x <0 \qquad \cite[\textbf{8.211(1)}]{grad}.
 \end{equation}
 For $x>0$, $\Ei$ is defined by \emph{the Cauchy principal value} (P.V.)  method:
\begin{equation*}
\Ei(x)= - \lim_{\epsilon \to 0}\left[ \int_{- x}^{-\epsilon}\frac{e^{-t}}{t}dt + \int_{\epsilon}^\infty
\frac{e^{-t}}{t}dt\right],  \quad  x>0   \qquad \cite[\textbf{8.211(2)}]{grad}.
\end{equation*}
Other useful representation is
$$\Ei(x)= \gamma + \log |x| + \sum_{k= 1}^\infty \frac{x^k}{k\,.\,k!}, \quad x\neq 0   \qquad \cite[\textbf{8.214(1,2)}]{grad}.$$
From above, we get:
\begin{equation}\label{eilog}
\frac{d}{dx}\Ei(\pm x)=\frac{e^{\pm x}}{x}\quad \mbox{and}\quad \lim_{x\to 0^+}\,\Ei(\pm x)-\log |x|=\gamma
= \int_0^1 \frac{1-e^{-t}}{t} dt-\int_1^\infty \frac{e^{-t}}{t} dt.
\end{equation}
We also have
\begin{equation}\label{eilog1}
\frac{1}{2}\,\left[e^{-x}\Ei(x)+ e^x \Ei(-x)\right]=-\int_0^\infty  \frac{t \cos(xt)}{1+t^2}\,dt, \quad x\in \Rset  \qquad \cite[\textbf{8.217(2)}]{grad}
\end{equation}
\item  The identity
\begin{equation}\label{latg}
i \log \frac{1-ie^x}{1+ie^x}=2 \arctan(e^x),  \quad x \in \Rset,
\end{equation}
is easily checkable by differentiation.

\item For Proposition \ref{pro4},  we need another two special functions, namely, the  {\it Riemann's zeta}  function, $\zeta$,  and the {\it polylogarithm functions} $\Li_n$, defined as follows:
\begin{eqnarray}
\zeta(s,a)&:=&\sum_{k=0}^\infty \frac{1}{(k+a)^s}, \quad \Re(s) > 1,\; - a \not\in \mathbb{N}; \label{rz}\\
\Li_s(z)&:=&\sum_{k=1}^\infty \frac{z^k}{k^s}, \quad |z|<1 , \; s \in \Cset , \label{pl}
\end{eqnarray}
with the properties, $\Li_s(1)=\zeta(s,1)$, and the representations
\begin{eqnarray}
\zeta(s,a) &=& \frac{1}{\Gamma(s)}\int_0^\infty e^{-a x} \frac{x^{s-1}}{1-e^{-x}} dx, \quad \Re(a)>0, \label{(8)}\\
\Li_s (z) &=& \frac{z}{\Gamma(s)}\int_1^\infty  \frac{\log^{s-1}(x)}{x(x-z)} dx, \quad z\in \Cset{\footnotesize \setminus} [1,\infty).\quad \cite[(7.188),\; p. \;236]{Lewin},\label{(81)}
\end{eqnarray}
where \eqref{(8)} is clearly a consequence of   \eqref{gag}. The name of the polylogarithm functions is due to the fact that $\Li_n,\, n=1,\,2, \dots ,$ can be defined recursively by
\begin{equation*}\label{li13}
\Li_n (z)=\int_0^z\frac{\Li_{n-1}(y)}{y}\,dy, \quad \Li_0(z)=\frac{z}{1-z}, \quad \Li_1(z)=-\log(1-z), \quad z<1.
\end{equation*}
By \eqref{latg} and \cite[Formulae (2.3) p. 38 and (4.29) p. 106]{Lewin}, we have
\begin{equation}\label{lig}
\Li_2(ie^s)-\Li_2(-ie^s) + s \log \frac{1-ie^s}{1+ie^s}= -2i \int_0^{s}\frac{x}{\cosh(x)}dx ,\quad s\in \Rset,
\end{equation}
where  $K=(\Li_2(i)-\Li_2(-i))/2i$  stands for the {\it Catalan constant}  $\approx$ 0.9159. In particular, we have:
$$   \Li_2(-i e^s)+ \Li_2(i e^s)  = 2i\sum_{k=1}^\infty (-1)^k \frac{e^{(2k-1)s}}{(2k-1)^2}. $$
\end{enumerate}

\end{document}